\documentstyle[12pt,a4,eepic,epsfig]{article}
\def\epsfile{\epsfig}

\title{
Algebra of screening operators 
\\
for the deformed $W_n$ algebra
}

\author{
Boris Feigin\thanks{L.D.Landau Institute for Theoretical Sciences,
			    Chernogolovka, 142432, Russia.},
Michio Jimbo\thanks{Department of Mathematics, Faculty of Science,
                            Kyoto University, Kyoto 606, Japan.},
Tetsuji Miwa\thanks{Institut Henri Poincar\'e and Ecole Normale
Superieure, on leave from Research Institute for Mathematical Sciences, 
                            Kyoto University, Kyoto 606, Japan.},
\\
Alexandr {Odesskii}$^*$,
and Yaroslav Pugai$^*$}

\date{February 25, 1997}      

\begin{document}

\renewcommand{\theequation}{\thesection.\arabic{equation}}

\newcommand{\End}{{\rm End}}
\newcommand{\Res}{{\rm Res}\,}
\newcommand{\id}{{\rm id}}
\newcommand{\sgn}{{\rm sgn}}
\newcommand{\wt}{{\rm wt}}
\newcommand{\nn}{\nonumber}
\newcommand{\omb}{\underline{\omega}}
\newcommand{\nb}{\underline{n}}
\newcommand{\eqref}[1]{(\ref{#1})}
\newcommand{\refeq}[1]{(\ref{eqn:#1})}
\newcommand{\be}{\begin{equation}}
\newcommand{\en}{\end{equation}}
\newcommand{\bea}{\begin{eqnarray}}
\newcommand{\ena}{\end{eqnarray}}
\newcommand{\bean}{\begin{eqnarray*}}
\newcommand{\enan}{\end{eqnarray*}}
\newcommand{\deq}{\stackrel{\rm def}{=}}
\newcommand{\lb}[1]{\label{eqn:#1}}
\newcommand{\mod}{~\hbox{mod}~}

\newcommand{\Bbb}{\bf}
\newcommand{\Z}{{\Bbb Z}}
\newcommand{\R}{{\Bbb R}}
\newcommand{\C}{{\Bbb C}}
\newcommand{\F}{{\cal F}}
\renewcommand{\S}{{\cal S}}
\newcommand{\ve}{\varepsilon}
\newcommand{\bve}{\bar{\varepsilon}}
\newcommand\e{\epsilon}
\newcommand\barep{\bar\ve}
\renewcommand\o{\omega}

\def\dlarrow{{\rm h}}
\def\duarrow{{\rm v}}
\newcommand{\lft}{\dlarrow}
\newcommand{\up}{\duarrow}
\newcommand{\br}[1]{{\langle #1 \rangle}}
\newcommand{\ket}[1]{{| #1 \rangle}}
\newcommand{\phit}{\tilde\phi}
\newcommand{\phib}{\bar\phi}
\newcommand{\qed}{\hfill \fbox{}\medskip}
\newcommand{\proof}{\medskip\noindent{\it Proof.}\quad }
\newcommand{\PP}{{\cal P}}
\def\iv#1{{d#1\over2\pi i#1}}
\def\hp{{\hat\pi}}
\def\deg{{\rm deg}}
\def\gl{{\lower.5ex\hbox{$>$}\atop\raise.5ex\hbox{$<$}}}
\def\lg{{\lower.5ex\hbox{$<$}\atop\raise.5ex\hbox{$>$}}}
\def\bX{{\overline X}}
\def\half{{1\over2}}
\def\frac#1#2{{#1\over#2}}
\def\ket#1{|#1\rangle}
\def\rst{{z^{(f_p)}_{k+p}=x^{s-p}w\phantom{m}(0\le p\le s)}}
\def\zvar{{z^{(1)}_i,\ldots,z^{(a)}_{i+m}}}
\def\prd#1{\prod_{1\le b#1 c\le a\atop i\le j\le i+m}}
\def\z#1#2{z^{(#1)}_{#2}}
\def\ab#1{#1^{(f_0)}_k\Bigg|_s}
\newtheorem{thm}{Theorem}[section]
\newtheorem{prop}[thm]{Proposition}
\newtheorem{lem}[thm]{Lemma}
\newtheorem{dfn}[thm]{Definition}

\maketitle

\begin{abstract}
We construct a family of intertwining operators (screening operators) 
between various Fock space modules over the deformed $W_n$ algebra. 
They are given as integrals involving 
a product of screening currents and elliptic theta functions. 
We derive a set of quadratic relations among the screening operators, 
and use them to construct a Felder-type complex 
in the case of the deformed $W_3$ algebra. 
\end{abstract}

\setcounter{section}{0}
\setcounter{equation}{0}

\section{Introduction}

The method of bosonization is known to be the 
most effective way of calculating the conformal blocks in 
conformal field theory. 
The basic idea in this approach is to 
realize the commutation relations for 
the symmetry algebra (such as the Virasoro or affine Lie algebras)  
and the chiral primary fields
in terms of operators acting on some bosonic Fock spaces. 
Quite often, the physical Hilbert space of the theory 
is not the total Fock space itself, but only a subquotient of it. 
In this case, it is necessary to `project out' the physical space 
from the Fock spaces by a cohomological method. 
In the case of the Virasoro minimal models, Felder \cite{Fel89}
introduced a two-sided complex 
\begin{equation}\label{eqn:Feld}
\cdots {\buildrel d \over \longrightarrow} \F^{(-1)}
{\buildrel d \over \longrightarrow} \F^{(0)}
{\buildrel d \over \longrightarrow} \F^{(1)}
{\buildrel d \over \longrightarrow} \F^{(2)}
{\buildrel d \over \longrightarrow}
\cdots,
\end{equation}
consisting of Fock spaces $\F^{(i)}$. 
As it turns out, the cohomology of this complex vanishes except 
at the $0$-th degree, and the remaining non-trivial cohomology 
affords the irreducible representation of the Virasoro algebra. 
The primary fields realized on the Fock spaces commute with 
the coboundary operator $d$, and hence make sense as operators on the 
cohomology space. 
Similar resolutions have been described for representations of
affine Lie algebras \cite{BeFe90,FeFr90,BMP90}. 
For an extensive review on this subject, 
the reader is referred to \cite{BMPrev}.  

It has been recognized that the 
idea of bosonization is quite fruitful also in 
non-critical lattice models \cite{collin,JM} 
and massive field theory \cite{Luk95}. 
The present work is motivated by the recent progress 
along this line, on the restricted solid-on-solid models \cite{LukPug2,AJMP}.
In \cite{LukPug2}, the Andrews-Baxter-Forrester (ABF) model was studied. 
Here the counterpart of the conformal primary fields are 
the vertex operators (half transfer matrices) 
which appear in the corner transfer matrix method. 
Just as in conformal field theory, the bosonization 
discussed in \cite{LukPug2} consists in two steps. 
The first step is to introduce a family of 
bosonic Fock spaces and realize the commutation relations of the 
vertex operators in terms of bosons. 
The second step is to realize 
the physical space of states of the model as the 
$0$-th cohomology of a complex of the type \eqref{eqn:Feld}. 
In fact, the analogy with conformal field theory goes further. 
Each Fock space has the structure of a module over 
the deformed Virasoro algebra (DVA) discovered in \cite{qVir}, 
where the deformation parameter $x$ ($0<x<1$) is the one which enters
the Boltzmann weights of the models. 
As was shown in \cite{LukPug2,JLMP}, the above complex 
is actually that of DVA modules, 
i.e., the operator $d$ commutes with the action of DVA. 
Felder's complex \eqref{eqn:Feld} is recovered in the limit 
$x\rightarrow 1$ (we shall refer to this as the conformal limit). 

The ABF models have $sl_n$ generalizations, the former being the case $n=2$. 
In the work \cite{AJMP}, 
the first step of the bosonization was carried over to the case of general $n$. 
However the second step was not addressed there. 
The aim of the present paper is to construct an analog of 
the complex \eqref{eqn:Feld} in the case $n=3$. 
In this situation the role of DVA is played by 
the deformed $W_3$ algebra introduced in \cite{qWN,FeFr95}. 
We shall also construct for general $n$ 
a family of intertwiners of deformed $W_n$ algebras (DWA), 
which we expect to be sufficient to construct the complex in the general case. 

In the conformal case, such a Felder-type complex for 
$\widehat{sl}_3$ was constructed in \cite{BMP90}. 
Strictly speaking, \cite{BMP90} discusses 
representations of affine Lie algebras, while our case corresponds
(in the limit) to those of $W_n$ algebras. 
In other words we are dealing with a coset theory rather than a 
Wess-Zumino-Witten theory. 
However the construction of the complex is practically the same for 
both cases. 

In the case $n\ge 3$, each component $\F^{(i)}$ of the complex is itself 
a direct sum of an infinite number of Fock spaces. 
The coboundary operator $d$ can be viewed as a collection of 
maps between various Fock spaces. 
We call these maps the screening operators. 
They are given in the form of an integral of a product of screening currents, 
multiplied by a certain kernel function expressed in terms of elliptic 
theta functions. 
The main result of this paper is the explicit construction of these 
screening operators. 
In comparison with the conformal case, 
a simplifying feature is that the screening operators 
can be expressed as products of more basic, mutually commuting operators. 
In the conformal case, such a multiplicative structure exists only 
`inside the contour integral' (see \cite{BMP90} and section 6
below). 
It has been pointed out \cite{FeFr95} 
that the screening currents satisfy the 
commutation relations of the elliptic algebra studied in \cite{FeOd95}. 
This connection turns out to be quite helpful in finding the 
basic operators referred to above and their commutation relations. 

Let us mention some questions that remain open. 
In order to ensure the nilpotency property $d^2=0$, 
the signs of the screening operators have to be chosen carefully. 
We have verified that this is possible for $n=3$. 
In the general case there are additional complications which we have not 
settled yet. 
More importantly, in this paper we do not discuss 
the cohomology of the complex, 
though we expect the same result persists as in the conformal case.
The construction of the complex in \cite{BMP90} is based on the 
one-to-one correspondence between intertwiners of Fock space
modules over $\widehat{sl}_n$, and the singular vectors in 
the Verma modules of $U_q(sl_n)$ with $q$ a root of unity. 
It would be interesting to search for 
an analog of the latter in the deformed situation.

The outline of this paper is as follows.
In section 2 we prepare the notation and the setting. 
Also the form of the complex to be constructed is briefly explained.
The construction of the screening operators for general $n$ is 
rather technical. 
To ease the reading, we first discuss in section 3 the case $n=3$ in detail. 
In section 4 we introduce the screening operators in general, 
and state their commutation relations. 
In section 5 we show that they commute with the action of DWA. 
In section 6 we briefly discuss the CFT limit
of the basic operators.
The text is followed by 3 appendices. 
In appendix A we discuss the condition 
when we construct intertwiners between two Fock modules. 
In appendix B we list the commutation relations for the basic operators 
that will be used to derive the quadratic relations of the screening 
operators. 
In appendix C we outline the proof that the screening operators commute with 
DWA. 

\setcounter{section}{1}
\setcounter{equation}{0}

\section{Preliminaries}

In this section we prepare the notation to be used in the text. 
Throughout this paper, we fix a positive integer $r\ge n+2$ and 
a real number $x$ with $0<x<1$. 
	
\subsection{Lie algebra $sl_n$}
Let us fix the notation concerning the Lie algebra $sl_n$.  

Let $\ve_i$ $(1\leq i \leq n)$ be an orthonormal basis in $\R^n$ 
relative to the standard inner product $(~,~)$. 
We set $\barep_i=\ve_i-\ve$,$\ve=(1/n)\sum_{j=1}^n\ve_j$. 
We shall denote by:
\begin{itemize}
\item $\o_i=\sum_{j=1}^i\barep_j$ the fundamental weights, 
\item $\alpha_i=\ve_i-\ve_{i+1}$ the simple roots, 
\item $\theta=\sum_{i=1}^{n-1}\alpha_i$ the maximal root,
\item $P=\sum_{i=1}^n\Z\barep_i$ the weight lattice, 
\item $Q=\sum_{i=1}^{n-1}\Z\alpha_i$ the root lattice,
\item $\Delta_+=\{\ve_i-\ve_j \mid 1\le i < j \le n-1 \}$ 
the set of positive roots, 
\item $\overline W\simeq S_n$ the classical Weyl group, 
\item $W\simeq \overline{W} |\!\!\!\times  Q$ the affine Weyl group. 
\end{itemize}
For an element $\gamma\in Q$, we set 
\[
|\gamma|=\sum_{i=1}^{n-1}c_i
\qquad \hbox{ for   }\gamma=\sum_{i=1}^{n-1}c_i\alpha_i.
\]
For a root $\alpha=\ve_i-\ve_j$, $r_\alpha$ signifies  
the corresponding reflection (often identified with the transposition 
$(ij)\in S_n$). 
We also write $s_i=r_{\alpha_i}$. 

\subsection{Bosons}
We recall from \cite{AJMP} our convention about the bosons. 
Let $\beta^j_m$ be the oscillators 
($1\le j\le n-1, m\in\Z\backslash\{0\}$) with the commutation relations 
\begin{eqnarray}
[\beta^j_m,\beta^k_{m'}]
&=&
m\frac{[(n-1)m]_x}{[nm]_x}\frac{[(r-1)m]_x}{[rm]_x}\delta_{m+m',0}
\qquad (j=k),
\label{eqn:2.1}\\
&=&
-mx^{sgn(j-k)nm}\frac{[m]_x}{[nm]_x}\frac{[(r-1)m]_x}{[rm]_x}\delta_{m+m',0}
\qquad (j\neq k).
\nonumber\\
&&\label{eqn:2.2}
\end{eqnarray}
Here the symbol $[a]_x$ stands for $(x^a-x^{-a})/(x-x^{-1})$. 
Define $\beta^n_m$ by
\begin{equation}
\sum_{j=1}^nx^{-2jm}\beta^j_m=0.
\label{eqn:2.3}
\end{equation}
Then the commutation relations \eqref{eqn:2.1},\eqref{eqn:2.2} are valid for 
all $1\le j,k\le n$. 

We also introduce the zero mode operators $P_\lambda,Q_\lambda$ indexed by 
$\lambda\in P$. 
By definition they are $\Z$-linear in $\lambda$ and satisfy 
\[
[iP_\lambda,Q_\mu]=(\lambda,\mu)
\qquad (\lambda,\mu\in P).
\]

We shall deal with the bosonic Fock spaces 
$\F_{l,k}$ ($l,k\in P$) generated by $\beta^j_{-m}$ ($m>0$) over the 
vacuum vectors $\ket{l,k}$:
\[
\F_{l,k}=\C[\{\beta^j_{-1},\beta^j_{-2},\cdots\}_{1\le j\le n}]\ket{l,k},
\]
where
\begin{eqnarray*}
&&\beta^j_m\ket{l,k}=0,\quad (m>0),
\\
&&P_\alpha\ket{l,k}=(\alpha,\sqrt{\frac{r}{r-1}}l-\sqrt{\frac{r-1}{r}}k)
\ket{l,k},
\\
&&\ket{l,k}=e^{i\sqrt{\frac{r}{r-1}}Q_l-i\sqrt{\frac{r-1}{r}}Q_k}\ket{0,0}.
\end{eqnarray*}
In what follows we set 
\[
\hat{\pi}_i=\sqrt{r(r-1)}P_{\alpha_i}.
\]
It acts on $\F_{l,k}$ as an integer,
\[
\hat{\pi}_i\bigl|_{\F_{l,k}}=(\alpha_i, rl-(r-1)k).
\]
In this paper we work on $\F_\lambda{\buildrel\rm def\over=}
\F_{l,\lambda}$ $(\lambda\in P)$
with a fixed value of $l\in P$.

\subsection{Screening currents}
We define the screening currents $\xi_j(u)$ ($j=1,\cdots,n-1$) by 
\begin{eqnarray}
\xi_j(u) \equiv U_j(z)=e^{i\sqrt{\frac{r-1}{r}}Q_{\alpha_j}}
z^{\frac{1}{r}\hat{\pi}_j+\frac{r-1}{r}}
:e^{\sum_{m\neq 0}\frac{1}{m}(\beta^j_m-\beta^{j+1}_m)(x^jz)^{-m}}:,
\label{eqn:2.4}
\end{eqnarray}
where the variable $u$ is related to $z$ via $z=x^{2u}$.

We shall need the following commutation relations between them. 
\begin{eqnarray}
&&\xi_j(u)\xi_j(v)=\frac{[u-v-1]}{[u-v+1]}\xi_j(v)\xi_j(u),
\label{eqn:co1}\\
&&\xi_j(u)\xi_{j\pm1}(v)=-\frac{[u-v+\frac{1}{2}]}{[u-v-\frac{1}{2}]}
\xi_{j\pm 1}(v)\xi_j(u),
\label{eqn:co2}\\
&&\xi_i(u)\xi_j(v)=\xi_j(v)\xi_i(u) \qquad (|i-j|>1).
\label{eqn:co3}
\end{eqnarray}
Here the symbol $[u]$ stands for the theta function 
satisfying
\bea
&&[u+r]=-[u]=[-u],\nonumber\\
&&[u+\tau]=-e^{2\pi i(u+{\tau\over2})/r}[u]
\hbox{ where $\tau={\pi i\over\log x}$}.\nonumber
\ena
Explicitly it is given by
\begin{eqnarray}
&&[u]=x^{u^2/r-u}\Theta_{x^{2r}}(x^{2u}),\\
&&\Theta_q(z)=(z;q)_\infty (qz^{-1};q)_\infty (q;q)_\infty, \\
&&(z;q)_\infty=\prod_{i=0}^\infty (1-zq^i).
\end{eqnarray}

Quite generally we say that an operator $X$ has weight $\nu$ 
if $X\F_\lambda\subset \F_{\lambda+\nu}$ for any $\lambda$. 
Then $\xi_j(u)$ has weight $-\alpha_j$. 
This implies 
\begin{equation}
\hat{\pi}_i\xi_j(u)
=\xi_j(u)\left(\hat{\pi}_i-(\alpha_i,\alpha_j)(1-r)\right).\label{ZMCOM}
\end{equation}

\subsection{The complex}\label{SEC3}

In this section we describe the form of the complex 
we are going to construct.  

Fix an integral weight $\Lambda\in P$ satisfying 
\bea
&&(\Lambda,\alpha_i)>0 \quad (i=1,\cdots,n-1),
\quad (\Lambda,\theta)<r. \label{3.0n}
\ena
Note that
\bea
0<(\Lambda,\alpha)<r\label{3.00}
\ena
for any positive root $\alpha$.
Consider the orbit of $\Lambda$ under the action of
the affine Weyl group $W$. 
An element of $W\Lambda$ can be written uniquely as 
\begin{equation}\label{eqn:aff}
\lambda=t_\gamma\sigma\Lambda=\sigma\Lambda+r\gamma,
\end{equation}
where $\sigma\in\overline{W}$ and $\gamma\in Q$. 
We assign a degree $\deg(\lambda)\in\Z$ to (\ref{eqn:aff}) by setting
\begin{equation}\label{eqn:deg}
\deg(\lambda)=l(\sigma)-2|\gamma|.
\end{equation}
Here $l(\sigma)$ denotes the length of $\sigma\in\overline{W}$.
(The right hand side of (\ref{eqn:deg}) is known as 
the modified length of $w=t_\gamma\sigma\in W$, see e.g.\cite{FeFr90,BMP90}.)
We shall construct a complex of the form
\begin{equation}\label{BRST}
\cdots {\buildrel d \over \longrightarrow} \F^{(-1)}_\Lambda
{\buildrel d \over \longrightarrow} \F^{(0)}_\Lambda
{\buildrel d \over \longrightarrow} \F^{(1)}_\Lambda
{\buildrel d \over \longrightarrow} \F^{(2)}_\Lambda
{\buildrel d \over \longrightarrow}
\cdots
\end{equation}
where 
\begin{equation}
\F^{(i)}_\Lambda=\bigoplus_{\lambda\in W\Lambda \atop \deg(\lambda)=i}
\F_{\lambda} \qquad (i\in \Z).
\end{equation}
Except for $n=2$, $\F^{(i)}_\Lambda$  
is a direct sum of an infinite number of Fock spaces. 

The coboundary map $d:\F^{(i)}_\Lambda\rightarrow\F^{(i+1)}_\Lambda$ 
can be viewed as a collection of operators 
\begin{equation}\label{eqn:dll}
d_{\lambda',\lambda}:
\F^{(i)}_\lambda \rightarrow \F^{(i+1)}_{\lambda'}
\end{equation}
associated with each pair $\lambda,\lambda'\in W\Lambda$ 
satisfying $\deg(\lambda')=\deg(\lambda)+1$. 
We shall impose a restriction on the possible pair 
$\lambda,\lambda'$ as explained below.

For a positive root $\alpha$ and an element $\lambda\in W\Lambda$, 
we define an integer $m_\alpha(\lambda)$ by 
\bea
&&0<m_\alpha(\lambda)<r,\\
&&m_\alpha(\lambda)\equiv(\lambda,\alpha)\bmod r.
\ena
In other words, if $\lambda=t_\gamma\sigma\Lambda$, then we have
\bea
m_\alpha(\lambda)=\cases{(\sigma\Lambda,\alpha)
&if $(\sigma\Lambda,\alpha)>0$;\cr
(\sigma\Lambda,\alpha)+r
&if $(\sigma\Lambda,\alpha)<0$.\cr}
\label{MAL}
\ena
Set
\bea\label{LA}
\lambda^\alpha=\lambda-m_\alpha(\lambda)\alpha
=\cases{t_\gamma r_\alpha\sigma\Lambda
&if $(\sigma\Lambda,\alpha)>0$;\cr
t_{\gamma-\alpha}r_\alpha\sigma\Lambda
&if $(\sigma\Lambda,\alpha)<0$.\cr}
\ena
Clearly the weight $\lambda^\alpha$ also belongs to the orbit $W\Lambda$.

\begin{dfn}
We say that an ordered pair $(\lambda,\lambda')$ is admissible if 
the following hold for some positive root $\alpha$. 
\begin{eqnarray}
&\lambda'=\lambda^\alpha,\\
&\deg(\lambda^\alpha)=\deg(\lambda)+1.
\end{eqnarray}
\end{dfn}

We set $d_{\lambda',\lambda}=0$ if $\lambda,\lambda'$ is not admissible.
Otherwise, write $\lambda'=\lambda^\alpha$ and 
\begin{equation}\label{eqn:X}
d_{\lambda^\alpha,\lambda}= X_\alpha(\lambda)
:\F_{\lambda}\longrightarrow \F_{\lambda^\alpha}.
\end{equation}
The construction of the complex is equivalent to finding an operator 
\eqref{eqn:X} for each admissible pair $\lambda,\lambda^\alpha$, 
so that we have $d^2=0$. 
We shall refer to \eqref{eqn:X} as a screening operator. 
We also require that the screening operators commute with the DWA generators.
In practice, we find it convenient to construct the screening operators in the form
\[
X_\alpha(\lambda)=s_\alpha(\lambda)\bX_\alpha(\lambda)
\]
where $s_\alpha(\lambda)=\pm 1$ is a sign factor.
In Section 4 we give both $s_\alpha(\lambda)$ and
${\overline X}_\alpha(\lambda)$ so that we have $d^2=0$
for the case $n=3$. The general case is incomplete because we could not
find a proper choice of the signs $s_\alpha(\lambda)$.

The construction of the screening operators \eqref{eqn:X} is based on the 
screening currents \eqref{eqn:2.4}. 
Let us consider the case where $\alpha$ in (\ref{eqn:X}) is a simple root 
$\alpha_j$. 
It turns out that the pair $\lambda,\lambda^{\alpha_j}$ is 
admissible for any $\lambda\in W\Lambda$ (see Lemma \ref{DL}). 
In this case the operator (\ref{eqn:X}) can be found as follows:
\begin{eqnarray}
&\bX_{\alpha_j}(\lambda)=X_j^{a}
\qquad (a=m_{\alpha_j}(\lambda)),
\label{eqn:simp}\\
&X_j=
\displaystyle \oint\frac{dz}{2\pi i z}\xi_j(u)
\frac{[u+\frac{1}{2}-\hat{\pi}_j]}{[u-\frac{1}{2}]}
\qquad (z=x^{2u}).
\label{eqn:simp2}
\end{eqnarray}
Here the integration is taken over the contour $|z|=1$. 
Notice that the kernel function 
$F(u)=[u+\frac{1}{2}-\hat{\pi}_j]/[u-\frac{1}{2}]$ 
has the quasi-periodicity 
\bea
&&F(u+\tau)=e^{2\pi i(1-\hat{\pi}_j)/r}F(u)
\label{QPER}
\ena
which ensures that the integrand of \eqref{eqn:simp2}
is a single valued function in $z$. 

For $n=2$, \eqref{eqn:simp} exhausts the possible screening operators. 
For $n\ge 3$ we must also construct operators corresponding to 
non-simple roots. 
As we shall see, they are given by similar (but more complicated) 
integrals over products of the screening currents. 

\setcounter{section}{2}
\setcounter{equation}{0}
\section{Case $n=3$}

Before embarking upon the construction of the complex in general, 
let us first elaborate on the case $n=3$. 
Hopefully this will make clear the main points of the construction. 

The following figure shows 
the configuration of the weights in the orbit $W\Lambda$ for $n=3$.

\medskip
\begin{center}
\fbox{\epsfile{file=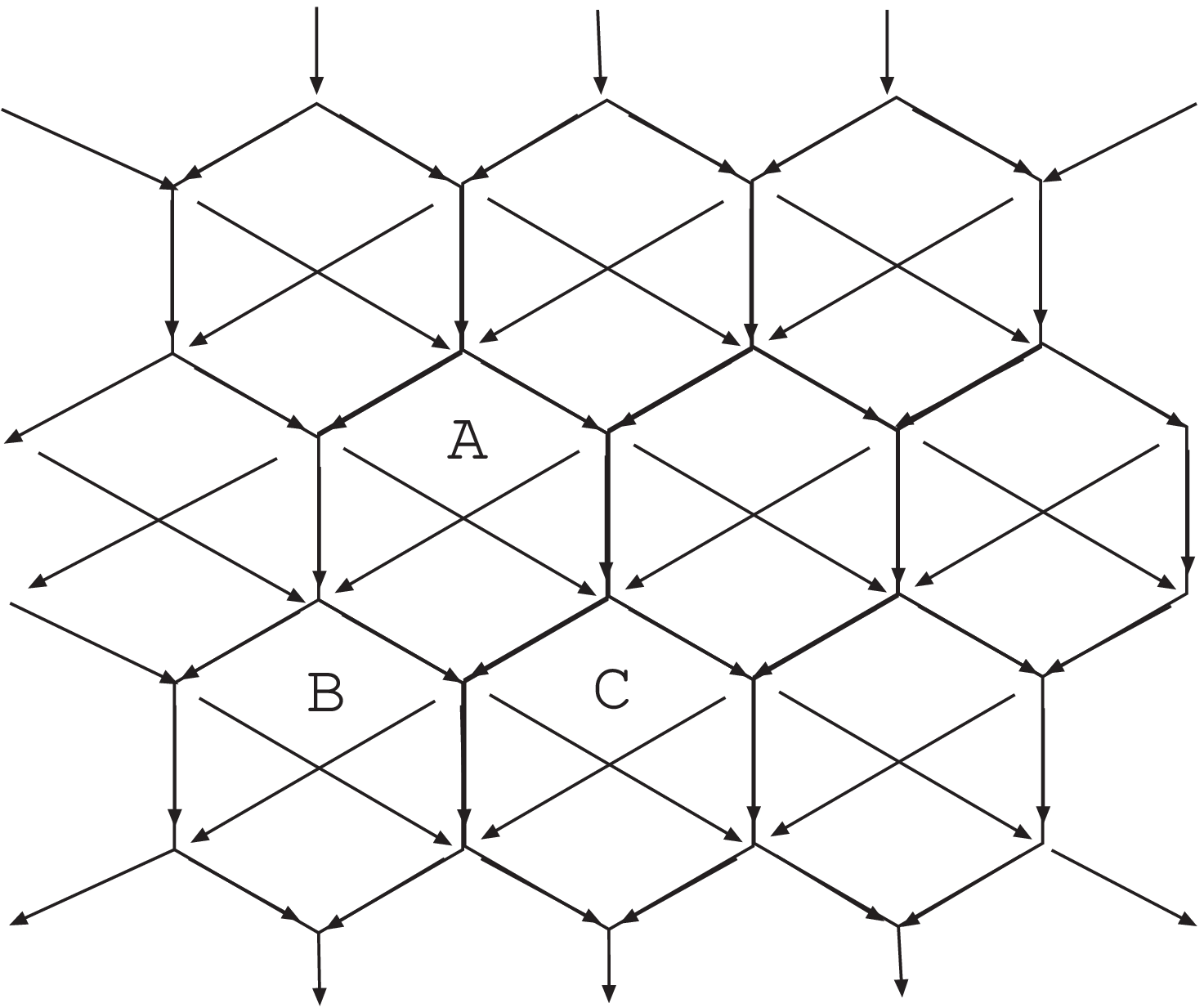,width=338pt,height=283pt}}
\end{center}
\medskip

\begin{quote}
Figure 1. \quad 
The orbit $W\Lambda$ for $n=3$.
It forms a hexagonal lattice, consisting of 
three types of basic hexagons $A,B,C$ and their translates
by $r$ times the root lattice $Q$.
\end{quote}

\bigskip

In the figure, each vertex represents a weight $\lambda\in W\Lambda$.
An arrow from $\lambda$ to $\lambda'$ indicates that 
the pair $(\lambda,\lambda')$ is admissible. 
As was mentioned before, $(\lambda,\lambda^\alpha)$ is 
always admissible for a simple root $\alpha=\alpha_1,\alpha_2$. 
For $n=3$ there is also the `third root' $\theta=\alpha_1+\alpha_2$. 
It turns out that $(\lambda,\lambda^\theta)$ is admissible if and only if
\[
\lambda=t_\gamma\sigma\Lambda
\quad\hbox{with }\sigma=s_1,s_2,s_1s_2s_1.
\]

The nilpotency $d^2=0$ leads to two types of relations
for the screening operators. 
The first type involves only 
screening operators corresponding to one simple root $\alpha_j$, and has the form
\[
X_j^a X_j^{r-a}=0.
\]
For $n=2$, this relation was proved in \cite{JLMP}. 
The same argument applies to show $X_j^r=0$ for any $j$. 

The second type involves the root $\theta=\alpha_1+\alpha_2$, and 
occurs for each square inside a hexagon (see Fig.1). 
Let us write the operator $\bX_{\theta}(\lambda)$ 
as $X_{12}^{(a)}$ ($a=m_\theta(\lambda)$), indicating that 
it has weight $-a\theta$. 
Let $a_i=(\Lambda,\alpha_i)$ ($i=1,2$), $a_0=r-a_1-a_2$, 
and suppose $\lambda=t_\gamma \sigma\Lambda$. 
If $\sigma=s_1$, then $m_\theta(\lambda)=a_2$, 
and the following relations must be satisfied:
\begin{eqnarray}
&&X_2^{r-a_1}X_1^{a_2}\pm X_{12}^{(a_2)}X_2^{a_0}=0,\label{eqn:d1}\\
&&X_2^{a_2}X_1^{r-a_1}\pm X_1^{a_0}X_{12}^{(a_2)}=0,\label{eqn:d2}\\
&&X_1^{r-a_0}X_2^{a_2}\pm X_{12}^{(a_2)}X_1^{a_1}=0,\label{eqn:d3}\\
&&X_1^{a_2}X_2^{r-a_0}\pm X_2^{a_1}X_{12}^{(a_2)}=0.\label{eqn:d4}
\end{eqnarray}
The operator $X^{(a)}_{12}$ and the signs $\pm$ in the commutation relations
will be given later.
The relations in the other cases $\sigma=s_2,s_1s_2s_1$ are obtained by 
permuting the upper indices $(a_0,a_1,a_2)\rightarrow (a_2,a_0,a_1)$
successively. 

In the conformal case $x=1$ (and $n=3$), 
the `third' screening operator $X_{\theta}(\lambda)$ 
satisfying the relations \eqref{eqn:d1}-\eqref{eqn:d4} 
was found in the work \cite{BMP90}. 
The construction in \cite{BMP90} is based on the observation that 
the screening operators $X_j$  for the simple roots
satisfy the same Serre relations 
as do the Chevalley generators of the quantum group 
$U_q(sl_n)$, with $q$ being a root of unity 
($q^{2r}=1$ in the present notation). 
With the aid of these relations, the operator
$X_{\theta}(\lambda)$ 
was expressed as a (non-commutative) polynomial in $X_1$ and $X_2$. 
This method does not easily generalize to the deformed case, since
it appears that there is no analog of the Serre relations between $X_1$ and $X_2$. 
Nevertheless there exists a family of operators, 
in terms of which the third screening operator can be 
written in a simple factorized form. 

Consider an operator of weight $-\theta$ of the form 
\[
\displaystyle\oint\!\!\oint\frac{dz_1}{2\pi i z_1}
\frac{dz_2}{2\pi i z_2}\xi_1(u_1)\xi_2(u_2)F(u_1,u_2),
\]
with some function $F(u_1,u_2)$ which is periodic in 
$u_i$ with period $r$. 
The integration is taken over $|z_1|=|z_2|=1$. 
Recall (\ref{ZMCOM}) and (\ref{QPER}).
In order that the integrand be single valued in $z_i$, we demand that 
\[
F(u_1+\tau,u_2)=e^{-2\pi i \hat{\pi}_1/r}F(u_1,u_2),
\qquad
F(u_1,u_2+\tau)=e^{-2\pi i (\hat{\pi}_2-1)/r}F(u_1,u_2),
\]
where $\tau=\pi i/\log x$.
Assume further that 
$F(u_1,u_2)$ is holomorphic except for possible simple poles at
$u_i=1/2$ and $u_1-u_2+1/2=0$. 
(As for the last pole, we have taken into account the commutation relation 
\eqref{eqn:co2}.)
If we regard the $\hat{\pi}_i$'s as constants, then the 
space of functions satisfying these conditions is $3$ dimensional.	
It is straightforward to find a spanning set of such functions. 
This motivates us to introduce the following family of operators parameterized 
by $k$: 
\begin{eqnarray}
&X_{12}(k)=(-1)^k
\displaystyle\oint\!\!\oint\frac{dz_1}{2\pi i z_1}
\frac{dz_2}{2\pi i z_2}\xi_1(u_1)\xi_2(u_2)
\nonumber\\
&\quad \times 
\displaystyle\frac{[u_1+k+\frac{1}{2}-\hat{\pi}_1]}{[u_1-\frac{1}{2}]}
\frac{[u_2-k-\frac{1}{2}-\hat{\pi}_2]}{[u_2-\frac{1}{2}]}
\frac{[u_1-u_2-k-\frac{1}{2}]}{[u_1-u_2+\frac{1}{2}]}.
\label{eqn:X12}
\end{eqnarray}

\begin{prop}\label{prop:X12}
\begin{eqnarray}
&&[X_{12}(k),X_{12}(l)]=0 \hbox{ for any } k,l.
\label{eqn:a1}
\\
&&X_1X_2=X_{12}(-1),\qquad X_2X_1=X_{12}(0),
\label{eqn:a2}\\
&&X_1X_{12}(k)=X_{12}(k-1)X_1,
\label{eqn:a3}\\
&&X_2X_{12}(k-1)=X_{12}(k)X_2.
\label{eqn:a4}
\end{eqnarray}
\end{prop}
The proof of these statements will be given later in the context of general $n$. 
Notice the periodicity relation
\begin{equation}\label{eqn:c}
X_{12}(k+r)=(-1)^{r-1}X_{12}(k).
\end{equation}

Set
\begin{equation}\label{eqn:x12}
X^{(a)}_{12}(k)=\prod_{b=1}^aX_{12}(k-b+1). 
\end{equation}
Then, for any non-negative integers $a,b$, we have 
\begin{eqnarray}
&X_1^{a+b}X_2^b=X^{(b)}_{12}(-a-1)X_1^a,\label{eqn:b1}\\
&X_1^aX_2^{a+b}=X_2^bX^{(a)}_{12}(-b-1),\label{eqn:b2}\\
&X_2^aX_1^{a+b}=X_1^bX_{12}^{(a)}(a+b-1),\label{eqn:b3}\\
&X_2^{a+b}X_1^{b}=X_{12}^{(b)}(a+b-1)X_2^a.\label{eqn:b4}
\end{eqnarray}
As an example, let us verify the first relation. 
Consider first the case $b=1$. From \eqref{eqn:a2} and \eqref{eqn:a3} we have
\begin{eqnarray*}
X_1^{a+1}X_2
&=&X_1^aX_{12}(-1)\\
&=&X_1^{a-1}X_{12}(-2)X_1
\\
&=&\cdots\\
&=&X_{12}(-a-1)X_1^a.
\end{eqnarray*}
The case of general $b$ follows immediately from this and the definition 
\eqref{eqn:x12}.
The other relations are derived in a similar manner. 

Now set 
\[
\bX_{\alpha_1+\alpha_2}(\lambda)=X_{12}^{(a)}(k), 
\qquad a=m_{\alpha_1+\alpha_2}(\lambda),~k=(\lambda,\alpha_1)-1.
\]
Comparing \eqref{eqn:b1}-\eqref{eqn:b4} with 
\eqref{eqn:d1}-\eqref{eqn:d4} 
and taking \eqref{eqn:c} into account, 
we see that the desired relations are satisfied up to sign. 

It remains to settle the issue about the signs $s_\alpha(\lambda)$. 
From the definition, the screening operators have the periodicity 
\[
\bX_\alpha(\lambda+\beta)=\bX_\alpha(\lambda)
\qquad \hbox{if }\beta\in r \Z\alpha_1 +2r \Z\alpha_2.
\]
Thus the signs can also be chosen according to the same periodicity.
A direct verification shows that the following is one possible solution
for the $s_\alpha(\lambda)$.

\medskip
\begin{center}
\fbox{\epsfile{file=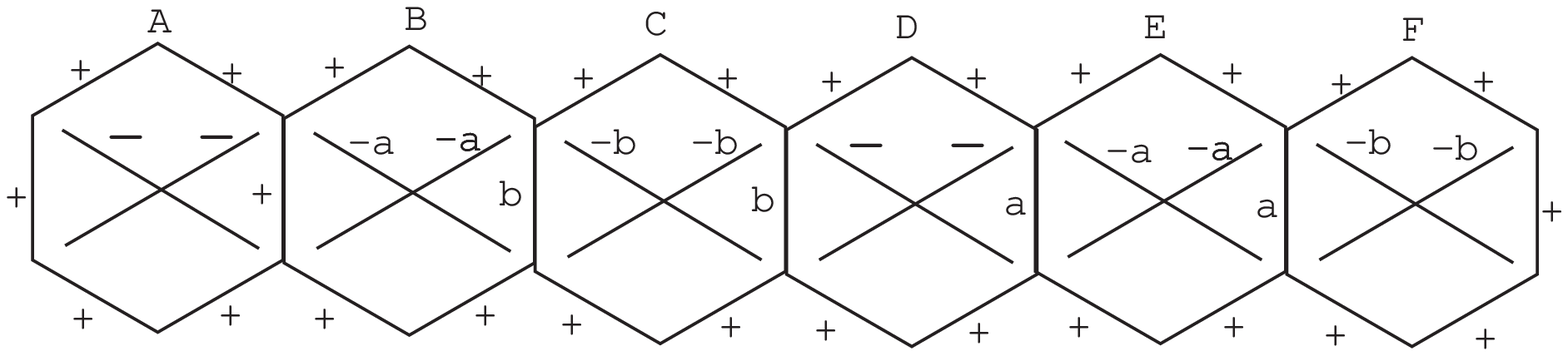,height=102pt,width=439pt}}
\end{center}
\medskip

\begin{quote}
Figure 2.\quad 
The choice of the signs $s_\alpha(\lambda)$.
We set $a=\ve_r^{k_1},b=\ve_r^{k_2},\ve_r=(-1)^{r-1}$.
The vertices at the top row correspond to the weights 
$A=\Lambda$, 
$B=r_1r_2\Lambda+\alpha_1$, 
$C=r_2r_1\Lambda+\alpha_1$, 
$D=\Lambda+\alpha_1-\alpha_2$,
$E=r_1r_2\Lambda+2\alpha_1-\alpha_2$, 
$F=r_2r_1\Lambda+2\alpha_1-\alpha_2$, 
with $r_j=r_{\alpha_j}$.
\end{quote}

\bigskip

\setcounter{section}{3}
\setcounter{equation}{0}
\section{Screening operators}

\subsection{Basic operators}
We are now in a position to introduce the operators which will 
play a basic 
role in the construction of screening operators for general $n$. 

Let $\alpha=\alpha_i+\cdots+\alpha_{i+m}$ be a positive root. 
We often write it as $\alpha_{i\cdots i+m}$. 
Define 
\bea\label{SO}
&&X_\alpha(k_1,\ldots,k_m)=\oint\cdots\oint\prod_{j=i}^{i+m}\iv{z_j}
{\xi_i(u_i)\over[u_i-\half]}\cdots{\xi_{i+m}(u_{i+m})\over[u_{i+m}-\half]}
\nonumber\\
&&\times\prod_{j=i+1}^{i+m}{[u_{j-1}-u_j]\over[u_{j-1}-u_j+\half]}
f_\alpha^{(k_1,\ldots,k_m)}(u_i,\ldots,u_{i+m};\hp_i,\ldots,\hp_{i+m}).
\ena
Here 
\bea
&&f^{(k_1,\ldots,k_m)}_\alpha(u_i,\ldots,u_{i+m};\hp_i,\ldots,\hp_{i+m})
=(-1)^{k_1+\cdots+k_m}
\nonumber\\
&&\prod_{l=1}^m
{[u_{i+l-1}-u_{i+l}-k_l-\half]\over[u_{i+l-1}-u_{i+l}]}
\prod_{l=0}^{m}[u_{i+l}-k_l+k_{l+1}-\half-\hp_{i+l}],
\label{SOS2}
\ena
and $k_0=-1,k_{m+1}=0$ is implied.  
The integrand of \eqref{SO} 
is a single-valued function in $z_j$ $(i\le j\le i+m)$.
To see this note that
\bea
\xi_i(u_i)\cdots\xi_{i+m}(u_{i+m})
\Bigl(\prod_{j=i}^{i+m}z_j^{-{1\over r}\hat\pi_j}\Bigr)
z_{i+m}^{-{r-1\over r}}
\ena
is single-valued.
When $\alpha=\alpha_j$ is a simple root, \eqref{SO} reduces to 
\eqref{eqn:simp}. 

The basic property of \eqref{SO} is the following commutativity.
\begin{thm}\label{TCOM}
For any $k_1,\cdots,k_m,p$ we have
\bea
[X_\alpha(k_1,\ldots,k_m),X_\alpha(k_1+p,\ldots,k_m+p)]=0.\label{COM}
\ena
\end{thm}

In view of this, we define for a non-negative integer $a$ 
\begin{dfn}\label{SO2}
\bea
X^{(a)}_\alpha(k_1,\ldots,k_m)&=&\prod_{b=1}^aX_\alpha(k_1-b+1,\ldots,k_m-b+1).
\label{MSC}
\ena
\end{dfn}
Sometimes, we abbreviate $X^{(a)}_{\alpha_{i\cdots i+m}}(k_1,\ldots,k_m)$
to $X^{(a)}_{i\cdots i+m}(k_1,\ldots,k_m)$.

\medskip

\noindent
{\it Proof of Theorem \ref{TCOM}.}\quad
Using the commutation relations \eqref{eqn:co1}-\eqref{eqn:co3} and 
\eqref{ZMCOM}, we can write the product 
$X_\alpha(k_1,\ldots,k_m)X_\alpha(k_1+p,\ldots,k_m+p)$
in the form 
\begin{eqnarray*}
&&\oint\!\!\cdots\!\!\oint\prod_{j=i}^{i+m}\iv{z_j}\iv{w_j}
\frac{\xi_i(u_i)}{[u_i-\half]}\frac{\xi_i(v_i)}{[v_i-\half]}
\cdots
\frac{\xi_{i+m}(u_{i+m})}{[u_{i+m}-\half]}
\frac{\xi_{i+m}(v_{i+m})}{[v_{i+m}-\half]}\\
&&\quad\times 
F(u_i,v_i,\cdots,u_{i+m},v_{i+m})
\qquad (z_k=x^{2u_k},w_k=x^{2v_k}), 
\end{eqnarray*}
and likewise for the product in the opposite order. 
Symmetrizing with respect to the integration variables and 
equating the integrand, we are led to prove the following equality. 
\bea
&&S\Biggl(
f_\alpha^{(k_1,\ldots,k_m)}
(u_i,\ldots,u_{i+m};\hp_i-(1-r),\hp_{i+1},\ldots,\hp_{i+m-1},
\nonumber\\
&&
\hp_{i+m}-(1-r))
f^{(k_1+p,\ldots,k_m+p)}_\alpha
(v_i,\ldots,v_{i+m};\hp_i,\ldots,\hp_{i+m})\nonumber\\
&&\times\prod_{j=i}^{i+m}{[u_j-v_j-1]\over[u_j-v_j]}
\prod_{j=i}^{i+m-1}{[u_j-v_{j+1}+\half]\over[u_j-v_{j+1}]}
\prod_{j=i+1}^{i+m}(-1){[u_j-v_{j-1}+\half]\over[u_j-v_{j-1}]}\Biggr)
\nonumber\\
&&=S\Biggl(
f^{(k_1+p,\ldots,k_m+p)}_\alpha(v_i,\ldots,v_{i+m};\hp_i-(1-r),\hp_{i+1},
\ldots,\hp_{i+m-1},
\nonumber\\
&&
\hp_{i+m}-(1-r))
f^{(k_1,\ldots,k_m)}_\alpha(u_i,\ldots,u_{i+m};\hp_i,\ldots,\hp_{i+m})
\nonumber\\
&&\times\prod_{j=i}^{i+m}{[u_j-v_j+1]\over[u_j-v_j]}
\prod_{j=i}^{i+m-1}(-1){[u_j-v_{j+1}-\half]\over[u_j-v_{j+1}]}
\prod_{j=i+1}^{i+m}{[u_j-v_{j-1}-\half]\over[u_j-v_{j-1}]}\Biggr).
\nonumber\\
\ena
Here the symbol $S$ means the symmetrization of $(u_j,v_j)$
for each $j=i,\ldots,i+m$.
This is equivalent to
\bea
&&A\Biggl(
\prod_{l=1}^m[u_{i+l-1}-u_{i+l}-k_l-\half]
[u_i+k_1+{3\over2}-\hp_i]
\biggl(\prod_{l=1}^{m-1}[u_{i+l}-k_l+k_{l+1}-\half-\hp_{i+l}]\biggr)
\nonumber\\
&&\times[u_{i+m}-k_m+\half-\hp_{i+m}]
\prod_{l=1}^m[v_{i+l-1}-v_{i+l}-k_l-p-\half]
[v_i+k_1+p+\half-\hp_i]
\nonumber\\&&
\biggl(\prod_{l=1}^{m-1}[v_{i+l}-k_l+k_{l+1}-\half-\hp_{i+l}]\biggr)
[v_{i+m}-k_m-p-\half-\hp_{i+m}]
\nonumber\\&&
\times\prod_{j=i}^{i+m}[u_j-v_j-1]
\prod_{j=i}^{i+m-1}[u_j-v_{j+1}+\half]
\prod_{j=i+1}^{i+m}[u_j-v_{j-1}+\half]\Biggr)
\nonumber\\&&
=A\Biggl(
\prod_{l=1}^m[u_{i+l-1}-u_{i+l}-k_l-\half]
[u_i+k_1+\half-\hp_i]
\biggl(\prod_{l=1}^{m-1}[u_{i+l}-k_l+k_{l+1}-\half-\hp_{i+l}]\biggr)
\nonumber\\&&
\times[u_{i+m}-k_m-\half-\hp_{i+m}]
\prod_{l=1}^m[v_{i+l-1}-v_{i+l}-k_l-p-\half]
[v_i+k_1+p+{3\over2}-\hp_i]
\nonumber\\&&
\times\biggl(\prod_{l=1}^{m-1}[v_{i+l}-k_l+k_{l+1}-\half-\hp_{i+l}]\biggr)
[v_{i+m}-k_m-p+\half-\hp_{i+m}]
\nonumber\\&&
\times\prod_{j=i}^{i+m}[u_j-v_j+1]
\prod_{j=i}^{i+m-1}[u_j-v_{j+1}-\half]
\prod_{j=i+1}^{i+m}[u_j-v_{j-1}-\half]\Biggr).
\nonumber\\
\ena
Here the symbol $A$ means the anti-symmetrization of
$(u_{i+j},v_{i+j})$ for each $j=0,\ldots,m$.

We prove this equality by induction on $m$.
The $\hbox{(LHS)}-\hbox{(RHS)}$ is a theta function 
of order $4$ in $u_i$.
Using the induction hypothesis, one can check that it
vanishes at $u_i=v_i,v_i\pm1$.
Taking into account the quasi-periodicity, we conclude that 
it must have a factor 
\bea
[u_i-v_i][u_i-v_i-1][u_i-v_i+1]
[u_i+2v_i-u_{i+1}-v_{i+1}+\half-\hp_i].
\ena
This is a contradiction unless $\hbox{(LHS)}-\hbox{(RHS)}=0$.
\qed

\subsection{Definition of screening operators}

Let us come to the definition of the screening operators. 
Let $(\lambda,\lambda^\alpha)$ be an admissible pair, 
with $\alpha=\alpha_{i\cdots i+m}$. 
We define $\bX_\alpha(\lambda):\F_\lambda\rightarrow\F_{\lambda^\alpha}$
by the formula 
\begin{equation}\label{SET1}
\bX_\alpha(\lambda)=X^{(a)}_\alpha(k_1,\ldots,k_m), 
\end{equation}
where
\bea
a&=&m_\alpha(\lambda),\label{SET2}\\
k_j&=&(\lambda,\alpha_{i\cdots i+j-1})-1.\label{SET3}
\ena
Note that on ${\cal F}_\lambda$ the operator
$\hp_j$ has the fixed value 
\bea
\hp_j\bigl|_{{\cal F}_\lambda}&=&(rl+(1-r)\lambda,\alpha_j)\nonumber\\
&=&(\lambda,\alpha_j)\bmod r.\label{SET4}
\ena

More explicitly the operator \eqref{SET1} is given as follows. 
\begin{prop}\label{prop:scr}
Notations being as in $(\ref{SOS2})$ and
$(\ref{SET1})$--$(\ref{SET4})$, we set 
\bea
&&f_\alpha^{(a)}(u^{(1)}_i,\ldots,u^{(a)}_i,
\ldots,u^{(1)}_{i+m},\ldots,u^{(a)}_{i+m})\nonumber\\
&&=S\Biggl(
\prod_{b=1}^a
f_\alpha^{(k_1-b+1,\ldots,k_m-b+1)}(u^{(b)}_i,\ldots,u^{(b)}_{i+m};
\hp_i-(a-b)(1-r),\hp_{i+1},\ldots,
\nonumber\\
&&
\ldots,\hp_{i+m-1},\hp_{i+m}-(a-b)(1-r))
\prod_{1\leq b<c\leq a\atop i\leq j\leq i+m}
{[u^{(b)}_j-u^{(c)}_j-1]\over[u^{(b)}_j-u^{(c)}_j]}
\nonumber\\
&&\times
\prod_{1\leq b<c\leq a\atop i\leq j\leq i+m-1}
{[u^{(b)}_j-u^{(c)}_{j+1}+{1\over2}]\over
[u^{(b)}_j-u^{(c)}_{j+1}]}\prod_{1\leq b<c\leq a\atop i+1\leq j\leq i+m}
(-1){[u^{(b)}_j-u^{(c)}_{j-1}+{1\over2}]
\over[u^{(b)}_j-u^{(c)}_{j-1}]}\Biggr)\nonumber\\
\label{LF1}\\
&&=S\Biggl(
\prod_{b=1}^a
f_\alpha^{(k_1-a+b,\ldots,k_m-a+b)}(u^{(b)}_i,\ldots,u^{(b)}_{i+m};
\hp_i-(a-b)(1-r),\hp_{i+1},\ldots,
\nonumber\\&&
\ldots,\hp_{i+m-1},\hp_{i+m}-(a-b)(1-r))
\prod_{1\leq b<c\leq a\atop i\leq j\leq i+m}
{[u^{(b)}_j-u^{(c)}_j-1]\over[u^{(b)}_j-u^{(c)}_j]}
\nonumber\\
&&\times
\prod_{1\leq b<c\leq a\atop i\leq j\leq i+m-1}
{[u^{(b)}_j-u^{(c)}_{j+1}+{1\over2}]\over
[u^{(b)}_j-u^{(c)}_{j+1}]}\prod_{1\leq b<c\leq a\atop i+1\leq j\leq i+m}
(-1){[u^{(b)}_j-u^{(c)}_{j-1}+{1\over2}]
\over[u^{(b)}_j-u^{(c)}_{j-1}]}\Biggr),\nonumber\\\label{LF2}
\ena
where $S$ means the symmetrization of $(u^{(1)}_j,\ldots,u^{(a)}_j)$
for each $j=i,\ldots,i+m$.
Then we have
\bea
&&\overline X_\alpha(\lambda)=
\oint\cdots\oint\prod_{1\leq b\leq a\atop i\leq j\leq i+m}
{dz^{(b)}_j\over2\pi iz^{(b)}_j}\nonumber\\
&&
\times \Biggl({\xi_i(u^{(1)}_i)\over[u^{(1)}_i-{1\over2}]}
\cdots{\xi_i(u^{(a)}_i)\over[u^{(a)}_i-{1\over2}]}\Biggr)
\cdots
\Biggl({\xi_{i+m}(u^{(1)}_{i+m})\over[u^{(1)}_{i+m}-{1\over2}]}
\cdots{\xi_{i+m}(u^{(a)}_{i+m})\over[u^{(a)}_{i+m}-{1\over2}]}\Biggr)
\nonumber\\
&&
\times\prod_{1\leq b<c\leq a\atop i\leq j\leq i+m}
\frac{[u^{(b)}_j-u^{(c)}_j]}{[u^{(b)}_j-u^{(c)}_j-1]}
\prod_{1\leq b,c\leq a\atop i\leq j\leq i+m-1}
\frac{[u^{(b)}_j-u^{(c)}_{j+1}]}
{[u^{(b)}_j-u^{(c)}_{j+1}+\frac{1}{2}]}
\nonumber\\
&&\times f^{(a)}_\alpha(u^{(1)}_i,\ldots,u^{(a)}_{i+m}).\label{eqn:scr}
\ena
The function $f^{(a)}_\alpha(u^{(1)}_i,\ldots,u^{(a)}_i,\ldots,
u^{(1)}_{i+m},\ldots,u^{(a)}_{i+m})$ has the following zeros:
\par\noindent $(i)$\quad
$u^{(b)}_j={1\over2}$ $(1\leq b\leq a,i\leq j\leq i+m)$
\par\noindent $(ii)$\quad
$u^{(b)}_j=u^{(c)}_{j+1}-\half=u^{(d)}_{j+1}+\half$
$(1\leq b,c,d\leq a,i\leq j\leq i+m-1)$
\par\noindent $(iii)$\quad
$u^{(b)}_j-\half=u^{(c)}_j+\half=u^{(d)}_{j+1}$
$(1\leq b,c,d\leq a,i\leq j\leq i+m-1)$
\par\noindent $(iv)$\quad
The function
\[
{f^{(a)}_\alpha(u^{(1)}_i,\ldots,u^{(a)}_{i+m})\over
\prod_{1\le b\le a\atop i\le j\le i+m}[u^{(b)}_j-\half]}
\]
is invariant under the simultaneous shift
$u^{(b)}_j\mapsto u^{(b)}_j+c$.
The function
\[
f^{(a)}_\alpha(u^{(1)}_i,\ldots,u^{(a)}_{i+m})
\prod_{1\le b,c\le a\atop i\le j\le i+m-1}[u^{(b)}_j-u^{(c)}_{j+1}]
\]
is holomorphic everywhere.
\end{prop}

\proof The equality of (\ref{LF1}) and  (\ref{LF2})
follows from the proof of Theorem \ref{TCOM}. 
Substituting \eqref{SET2}-\eqref{SET4} into \eqref{LF1}
and comparing with the definition \eqref{SOS2}, 
we see that the function 
$f_\alpha^{(k_1-b+1,\ldots,k_m-b+1)}$ in \eqref{LF1}
has a factor $\prod_{j=i+1}^{i+m}[u^{(b)}_j-1/2]$. 
Therefore $f^{(a)}_\alpha$ has a zero at
$u^{(b)}_j={1\over2}$ for $1\leq b\leq a$ and $i+1\leq j\leq i+m$. 
A similar argument applied to (\ref{LF2}) shows that it is true
also for $j=i$. 

For the proof of the rest of Proposition \ref{prop:scr},
we recall the definition of the elliptic algebra in \cite{FeOd95}.
(See also Appendix \ref{app:2}.)
The adaptation to our context goes as follows.

For $\hp_i\in\Z/r\Z$ $(1\le i\le n-1)$ and
$\gamma=\sum_{i=1}^{n-1}a_i\alpha_i$ $(a_i\in\Z_{\ge0})$
define $F^{(\hp_1,\ldots,\hp_{n-1})}_\gamma$
to be the set of functions in $|\gamma|=\sum_{i=1}^{n-1}a_i$
variables
\[
u^{(b)}_i\quad(1\le i\le n-1,1\le b\le a_i)
\]
satisfying the following properties:
$f(u^{(1)}_1,\ldots,u^{(a_{n-1})}_{n-1})$ is
\bea
&(P1)&\quad
\hbox{meromorphic in $u^{(b)}_i\in\C$,}\nonumber\\
&(P2)&\quad
\hbox{symmetric in each set of variables $(u^{(1)}_i,\ldots,u^{(a_i)}_i)$,}\nonumber\\
&(P3)&\quad
\hbox{quasi-periodic in each variable $u^{(b)}_i$
in the following sense:}\nonumber\\
&&f(u^{(b)}_i+r)=-f(u^{(b)}_i),\nonumber\\
&&f(u^{(b)}_i+\tau)=-e^{2\pi i(u^{(b)}_i+{\tau\over2}+
{(\gamma,\alpha_i)-1\over2}-\hp_i)/r}f(u^{(b)}_i),\nonumber\\
&(P4)&\quad
\hbox{holomorphic except for (at most) simple poles at
$u^{(b)}_i=u^{(c)}_{i+1}$,}\nonumber\\
&(P5)&\quad
\hbox{zero at one of the following:}\nonumber
\ena
\begin{itemize}
\item $u^{(b)}_i=u^{(c)}_{i+1}-\half=u^{(d)}_{i+1}+\half$
\item $u^{(b)}_{i+1}=u^{(c)}_i-\half=u^{(d)}_i+\half$
\end{itemize}

Note that if $a_i=0$ for some $i$ then $\hp_i$ does not appear in the
condition for $f$. If so, we sometimes abbreviate $\hp_i$ from the notation
$F^{(\hp_1,\ldots,\hp_{n-1})}_\gamma$.

Let
\[
\gamma^{(1)}=\sum_{i=1}^{n-1}a^{(1)}_i\alpha_i,
\qquad
\gamma^{(2)}=\sum_{i=1}^{n-1}a^{(2)}_i\alpha_i,
\]
and
let $\hp'_i$ be given by 
\bea\label{SHT}
\hp'_i=\hp_i-(\gamma^{(2)},\alpha_i)(1-r).
\ena
Then one can show \cite{FeOd95} that there exists an associative
mapping (the $*$-product)
\[
F^{(\hp'_1,\ldots,\hp'_{n-1})}_{\gamma^{(1)}}
\otimes
F^{(\hp_1,\ldots,\hp_{n-1})}_{\gamma^{(2)}}
\rightarrow F^{(\hp_1,\ldots,\hp_{n-1})}_{\gamma^{(1)}+\gamma^{(2)}},
\]
which sends $f\otimes g$ to $f*g$ where
\bea
&&(f*g)(u^{(1)}_1,\ldots,u^{(a^{(1)}_1)}_1,
v^{(1)}_1,\ldots,v^{(a^{(2)}_1)}_1,
\ldots,
u^{(1)}_{n-1},\ldots,u^{(a^{(1)}_{n-1})}_{n-1},
v^{(1)}_{n-1},\ldots,
\nonumber\\
&&\ldots, v^{(a^{(2)}_{n-1})}_{n-1})
=
S\Biggl(
f(u^{(1)}_1,\ldots,u^{(a^{(1)}_{n-1})}_{n-1})
g(v^{(1)}_1,\ldots,v^{(a^{(2)}_{n-1})}_{n-1})
\nonumber\\
&\times&
\prod_{
{{\scriptstyle 1\le i\le n-1}
\atop
{\scriptstyle 1\le b\le a^{(1)}_i}}
\atop
{\scriptstyle 1\le c\le a^{(2)}_i}}
{[u^{(b)}_i-v^{(c)}_i-1]\over[u^{(b)}_i-v^{(c)}_i]}
\prod_{
{{\scriptstyle 1\le i\le n-2}
\atop
{\scriptstyle 1\le b\le a^{(1)}_i}}
\atop
{\scriptstyle 1\le c\le a^{(2)}_{i+1}}}
{[u^{(b)}_i-v^{(c)}_{i+1}+\half]\over[u^{(b)}_i-v^{(c)}_{i+1}]}
\nonumber\\
&\times&
\prod_{
{{\scriptstyle 2\le i\le n-1}
\atop
{\scriptstyle 1\le b\le a^{(1)}_i}}
\atop
{\scriptstyle 1\le c\le a^{(2)}_{i-1}}}
(-1){[u^{(b)}_i-v^{(c)}_{i-1}+\half]\over[u^{(b)}_i-v^{(c)}_{i-1}]}
\Biggr).\nonumber
\ena
Here $S$ means  the symmetrization of $(u^{(1)}_i,\ldots,u^{(a^{(1)}_i)},
v^{(1)}_i,\ldots,v^{(a^{(2)}_i)}_i)$ for each $i=1,\ldots,n-1$. 
Because of the symmetrization the pole at $u^{(b)}_i=v^{(c)}_i$
is canceled. The property (P5) is preserved by the $*$-product.

Suppose that
$F(u^{(1)}_1,\ldots,u^{(a^{(1)}_{n-1})}_{n-1};\hp_1,\ldots,\hp_{n-1})$
is a function of the variables
$(u^{(1)}_1,\ldots,u^{(a^{(1)}_{n-1})}_{n-1})$
and $(\hp_1,\ldots,\hp_{n-1})$, and belongs to
$F^{(\hp_1,\ldots,\hp_{n-1})}_{\gamma^{(1)}}$
for any choice of $(\hp_1,\ldots,\hp_{n-1})\in\Z^{n-1}$.

Note that $(\hp_1,\ldots,\hp_{n-1})$ in the expression
$F^{(\hp_1,\ldots,\hp_{n-1})}_{\gamma^{(1)}}$
is considered as an element of $(\Z/r\Z)^{n-1}$.
Therefore, the factor $1-r$ in the shift (\ref{SHT}) could be simply
$1$. However, the function $F$ may have dependence on
$(\hp_1,\ldots,\hp_{n-1})\in\Z^{n-1}$. In fact, we use 
(\ref{SHT}) in this form in (\ref{FG}) below.

We define an operator $X(F)$ by
\bea
&&X(F)=
\oint\iv{z^{(1)}_1}\cdots\oint\iv{z^{(a^{(1)}_{n-1})}_{n-1}}
{\xi_1(u^{(1)}_1)\over[u^{(1)}_1-\half]}\cdots
{\xi_{n-1}(u^{(a^{(1)}_{n-1})}_{n-1})\over[u^{(a^{(1)}_{n-1})}_{n-1}-\half]}
\nonumber\\
&&\times\prod_{1\le i\le n-1\atop1\le b<c\le a^{(1)}_i}
{[u^{(b)}_i-u^{(c)}_i]\over[u^{(b)}_i-u^{(c)}_i-1]}
\prod_{{\scriptstyle1\le i\le n-1\atop
\scriptstyle1\le b\le a^{(1)}_i}\atop
\scriptstyle1\le c\le a^{(1)}_{i+1}}
{[u^{(b)}_i-u^{(c)}_{i+1}]\over[u^{(b)}_i-u^{(c)}_{i+1}+\half]}\nonumber\\
&&\times
F(u^{(1)}_1,\ldots,u^{(a^{(1)}_{n-1})}_{n-1};\hp_1,\ldots,\hp_{n-1}).
\ena
Similarly we define $X(G)$ from
$G(v^{(1)}_1,\ldots,v^{(a^{(2)}_{n-1})}_{n-1};\hp_1,\ldots,\hp_{n-1})$.

Set
\bea
f(u^{(1)}_1,\ldots,u^{(a^{(1)}_{n-1})}_{n-1})&=&
F(u^{(1)}_1,\ldots,u^{(a^{(1)}_{n-1})}_{n-1};\hp'_1,\ldots,\hp'_{n-1}),
\label{FG}\\
g(v^{(1)}_1,\ldots,v^{(a^{(2)}_{n-1})}_{n-1})&=&
G(v^{(1)}_a,\ldots,v^{(a^{(1)}_{n-1})}_{n-1};\hp_1,\ldots,\hp_{n-1}).
\nonumber\\
\ena
Then we have $X(F)X(G)=X(f*g)$.

{}From this follows (\ref{eqn:scr}). Then, (ii) and (iii) are nothing but
the condition (P5).

Let us prove (iv).
Assume that $f$ has zeros at $u^{(b)}_i=\half$ 
$(1\le i\le n-1,1\le b\le a^{(1)}_i)$.
Set
\bea
{\bar f}(u^{(1)}_1,\ldots,u^{(a^{(1)}_{n-1})}_{n-1})=
{f(u^{(1)}_a,\ldots,u^{(a^{(1)}_{n-1})})\over
\prod_{1\le i\le n-1 \atop 1\le b\le a^{(1)}_i}[u^{(b)}_i-\half]}.
\ena
From (P3) we have the following periodicity with respect to each
$u^{(b)}_i$
\bea
{\bar f}(u^{(b)}_i+r)&=&{\bar f}(u^{(b)}_i),\nonumber\\
{\bar f}(u^{(b)}_i+\tau)&=&e^{2\pi i({(\gamma,\alpha_i)\over2}-\hp_i)/r}
{\bar f}(u^{(b)}_i).\nonumber
\ena
Fix generic $(u^{(1)}_1,\ldots,u^{(a^{(1)}_{n-1})}_{n-1})$ and consider
a simultaneous shift of $\bar f$
\begin{equation}
{\bar f}_s(c)={\bar f}(u^{(1)}_1+c,\ldots,u^{(a^{(1)}_{n-1})}_{n-1}+c).
\end{equation}
If
\bea
{(\gamma,\gamma)\over2}\equiv
\sum_{i=1}^{n-1}a_i\hp_i\bmod r,\label{eq:C}
\ena
then ${\bar f}_s(c)$ is doubly periodic in $c$.
In Proposition \ref{prop:scr} we choose
$\gamma=(\lambda,\alpha)\alpha$,
and then (\ref{eq:C}) is valid.
From (P4) ${\bar f}_s$ has no pole in $c$, and therefore it is constant.
This completes the proof of Proposition \ref{prop:scr}.
\qed

Note that the statement (iv) of Proposition \ref{prop:scr} corresponds to the 
closeness of Felder's contour in the CFT limit.

\subsection{Quadratic relations}

With the screening operators introduced above, let us examine 
the nilpotency property $d^2=0$.
Unfortunately, we could not find a complete solution
to this problem in the general case $n\ge4$.
In Appendix \ref{app:1} we show the following. 
\begin{thm}\label{thm:comsq}
Suppose that $\alpha,\beta\in\Delta_+$ and 
\bea
(\lambda,\lambda^\alpha,\lambda^{\alpha,\beta})\label{C1}
\ena
is admissible, i.e., both $(\lambda,\lambda^\alpha)$ and
$(\lambda^\alpha,\lambda^{\alpha,\beta})$ are admissible.

$(i)$ If $(\alpha,\beta)=2$, then $\alpha=\beta$ and it is a simple root.

$(ii)$ Otherwise, there exists $\alpha'(\not=\alpha),\beta'\in\Delta_+$
such that
\bea
(\lambda,\lambda^{\alpha'},\lambda^{\alpha',\beta'})\label{C2}
\ena
is admissible and
\bea
\lambda^{\alpha,\beta}=\lambda^{\alpha',\beta'}.\label{C3}
\ena
The pair $(\alpha',\beta')$ is uniquely determined by this condition.
\end{thm}
In the second case, 
we say that the set of admissible weights
(\ref{C1}) and (\ref{C2}) satisfying (\ref{C3}) form 
a {\it commutative square}.

Consider the case (i). 
If $\alpha=\alpha_j$ is a simple root and $a=m_\alpha(\lambda)$, then 
we have $m_\alpha(\lambda^\alpha)=r-a$. 
Hence $d^2=0$ is ensured by the relation
\[
X_j^{r-a}X_j^a=X_j^r=0.
\]
This has been proved in \cite{JLMP}. 

In the case (ii), we need the following. 
\begin{thm}\label{thm:5.5}
For each commuting square, 
we have the identity of screening operators
\bea\label{ID}
\bX_\beta(\lambda^\alpha)\bX_\alpha(\lambda)=
s_\lambda(\alpha,\beta;\alpha'\beta')
\bX_{\beta'}(\lambda^{\alpha'})\bX_{\alpha'}(\lambda)
\ena
where
$s_\lambda(\alpha,\beta;\alpha'\beta')=\pm1$.
\end{thm}
The sign factor arises here because of the (anti-)periodicity 
property of the operator $X_\alpha(k_1,\ldots\ldots,k_m)$,  
\bea
X_\alpha(k_1,\ldots,k_m)|_{k_i\rightarrow k_i+r}
=\ve_r X_\alpha(k_1,\ldots,k_m),
\qquad (\ve_r=(-1)^{r+1}).
\ena
The precise formula for $s_\lambda(\alpha,\beta;\alpha'\beta')=\pm1$
will be given below. To have $d^2=0$ we must choose
the signs $s_\alpha(\lambda)$ appropriately. Theorem \ref{thm:5.5}
reduces the problem to finding $s_\alpha(\lambda)$ satisfying
\[
s_\beta(\lambda^\alpha)s_\alpha(\lambda)=-s_\lambda(\alpha,\beta;\alpha',
\beta')s_{\beta'}(\lambda^{\alpha'})s_{\alpha'}(\lambda).
\]
We have no solution except for the special cases $n=2,3$.

The assertion \eqref{ID} amounts to a number of identities 
of theta functions. 
These are derived in appendix \ref{app:2}. 
Below we shall indicate which identities are used in each case. 

\medskip
\noindent{\bf Case $(\alpha,\beta)=0$}\quad
In this case, $\alpha'=\beta$, $\beta'=\alpha$
hold (see Appendix \ref{app:1}). 
The assertion \eqref{ID} is nothing but (\ref{FR3}) and (\ref{FR4}). 
We can show that $s_\lambda(\alpha,\beta;\alpha',\beta')=1$. 

\medskip
\noindent{\bf Case $(\alpha,\beta)=\pm 1$}\quad
From the case-by-case analysis of Appendix \ref{app:1}, 
we see that altogether there are $8$ cases to consider. 
In the following we set
$\gamma_1=\alpha_{i\cdots i+l}$ and $\gamma_2=\alpha_{i+l+1\cdots i+l+m}$.

{\sl Case $A_+$}: $\alpha=\gamma_1,\beta=\gamma_2,
m_\alpha(\lambda)>m_\beta(\lambda^\alpha)$.

\smallskip
{\sl Case $B_+$}: $\alpha=\gamma_1,\beta=\gamma_2,
m_\alpha(\lambda)<m_\beta(\lambda^\alpha)$.

\smallskip
{\sl Case $C_+$}: $\alpha=\gamma_2,\beta=\gamma_1,
m_\alpha(\lambda)>m_\beta(\lambda^\alpha)$.

\smallskip
{\sl Case $D_+$}: $\alpha=\gamma_2,\beta=\gamma_1,
m_\alpha(\lambda)<m_\beta(\lambda^\alpha)$.

\smallskip
{\sl Case $A_-$}: $\alpha=\gamma_1+\gamma_2,\beta=\gamma_1$.

\smallskip
{\sl Case $B_-$}: $\alpha=\gamma_2,\beta=\gamma_1+\gamma_2$.

\smallskip
{\sl Case $C_-$}: $\alpha=\gamma_1+\gamma_2,\beta=\gamma_2$.

\smallskip
{\sl Case $D_-$}: $\alpha=\gamma_1,\beta=\gamma_1+\gamma_2$.

The cases $X_+$ and $X_-$ $(X=A,B,C,D)$ form a commutative square.
We will prove this statement and show the corresponding equality of
the screening operators. 

Set
\bea
\kappa_\alpha(\lambda)={m_\alpha(\lambda)-(\sigma\Lambda,\alpha)\over r}.
\ena
To get (\ref{ID}) for $\alpha=\beta'=\gamma_1$,
$\alpha'=\gamma_1+\gamma_2$, $\beta=\gamma_2$, we use (\ref{CS3}) with
\bea
a&=&m_{\alpha'}(\lambda),\nonumber\\
b&=&m_\alpha(\lambda)-m_{\alpha'}(\lambda),\nonumber\\
k_j&=&(\lambda,\alpha_{i\cdots i+j-1})-1\quad(1\leq j\leq l
\hbox{ or }l+2\leq j\leq l+m).\nonumber
\ena
The signature in (\ref{FR3}) is given by
\bea
s_\lambda(\alpha,\beta;\alpha',\beta')
=\ve_r^{\kappa_\alpha(\lambda)|\beta|m_{\alpha'}(\lambda)}.
\ena
 
To get (\ref{ID}) for $\alpha=\gamma_1$, $\beta=\alpha'=\gamma_2$,
$\beta'=\gamma_1+\gamma_2$, we use (\ref{CS4}) with
\bea
a&=&m_{\alpha'}(\lambda),\nonumber\\
b&=&m_\alpha(\lambda),\nonumber\\
k_j&=&(\lambda,\alpha_{i\cdots i+j-1})-1\quad(1\leq j\leq l
\hbox{ or }l+2\leq j\leq l+m).\nonumber
\ena
The signature in (\ref{FR3}) is given by
\bea
s_\lambda(\alpha,\beta;\alpha',\beta')
=\ve_r^{\kappa_\alpha(\lambda)|\alpha'|m_{\alpha}(\lambda)}.
\ena

To get (\ref{ID}) for $\alpha=\gamma_2$, $\beta=\alpha'=\gamma_1$,
$\beta'=\gamma_1+\gamma_2$, we use (\ref{CS1}) with
\bea
a&=&m_{\alpha'}(\lambda),\nonumber\\
b&=&m_\alpha(\lambda),\nonumber\\
k_j&=&(\lambda,\alpha_{i\cdots i+j-1})-1\quad(1\leq j\leq l),\nonumber\\
k_j&=&(\lambda,\alpha_{i+l+1\cdots i+j-1})-1\quad(l+2\leq j\leq l+m).\nonumber
\ena
The signature in (\ref{FR3}) is given by
\bea
s_\lambda(\alpha,\beta;\alpha',\beta')
=\ve_r^{\kappa_{\alpha'}(\lambda)|\alpha|m_{\alpha}(\lambda)}.
\ena

To get (\ref{ID}) for $\alpha=\beta'=\gamma_2$, $\beta=\gamma_1$,
$\alpha'=\gamma_1+\gamma_2$, we use (\ref{CS2}) with
\bea
a&=&m_{\alpha'}(\lambda),\nonumber\\
b&=&m_\alpha(\lambda)-m_{\alpha'}(\lambda),\nonumber\\
k_j&=&(\lambda,\alpha_{i\cdots i+j-1})-1\quad(1\leq j\leq l),\nonumber\\
k_j&=&(\lambda,\alpha_{i+l+1\cdots i+j-1})-1\quad(l+2\leq j\leq l+m).\nonumber
\ena
The signature in (\ref{FR3}) is given by
\bea
s_\lambda(\alpha,\beta;\alpha',\beta')
=\ve_r^{(\kappa_{\alpha}(\lambda)+\kappa_{\alpha'}(\lambda))
|\alpha|m_{\alpha'}(\lambda)}.
\ena

\setcounter{section}{4}
\setcounter{equation}{0}
\section{Operators $\bX_{\alpha}(\lambda)$ as intertwiners of DWA}

In this section we demonstrate that the screening operators 
introduced above are the intertwining operators for the DWA. 
A special case of the statement 
has been proved \cite{qWN,FeFr95} for the screening
operators acting in the vacuum module, where the 
theta function factor becomes unit and
the screening operator has a particularly simple form:
\bea
\bX_{\alpha_{j}}=\oint
{dz\over2\pi iz}\ U_j(z). 
\ena
The proof of \cite{qWN,FeFr95} was based on the fact that the
screening currents commute with DWA generators 
up to a total difference.
As it was explained \cite{JLMP}
on the example of the $sl_2$ case,
in dealing with general $(\alpha,\lambda)$
one needs to be sure that 
the additional theta function terms do not 
lead to nonvanishing contributions to the 
commutator of the DWA generators and $\bX_{\alpha}(\lambda)$.
We prove that this property is guaranteed by
the relations of the elliptic algebra of screening 
operators. 

\subsection{Deformed $W$ Algebra}
The DWA can be constructed as the subalgebra in the 
universal enveloping algebra of the Heisenberg algebra 
of operators $\beta_n$ and $P_\alpha$. Let us introduce 
the local bosonic fields \cite{qWN,FeFr95} 
\footnote{The difference between our notations
and those in \cite{qWN,FeFr95} is that we are 
working with the zero mode $\PP$ shifted as
$
\PP \longrightarrow \PP-\frac{1}
{\sqrt{r(r-1)}}\sum_{j=1}^{n-1}\vec{\omega}_j\ .
$
In the conformal limit $x\rightarrow 1$ this corresponds 
to the transform from the complex plane to the annulus.
The parameters $p,t,q$ in \cite{qWN} are related with 
$x$ and $r$
as $q=x^{-2r},\ t=x^{-2(r-1)},\ p=x^{-2}$.
}
\bea
\Lambda_j(z)=x^{-2\sqrt{r(r-1)}P_{\bve_j}}:\exp{(\sum_{m\neq 0}
\frac{(x^{rm}-x^{-rm})}{m}\beta^j_m z^{-m})
}:
\ena
The explicit realization of the generators of 
DWA in terms of these fields
is given by means of the non-linear transformation:
\bea
&&:(x^{-2z\partial_z}-\Lambda_1(z))(x^{-2z\partial_z}-\Lambda_2(zx^{2}))
\cdots (x^{-2z\partial_z}-\Lambda_n(zx^{2n-2})):
\nonumber\\
&&=\sum_{j=0}^{n}(-1)^j W^{(j)}(zx^{j-1})x^{-2(n-j)z\partial_z} \ ,
\label{DWAmiura}
\ena
where $W^{(0)}(z)\equiv 1$.
In the limit $x\rightarrow 1$ this formula 
leads to the quantum Miura transformation describing the 
bosonic realization of the $W_n$-algebra  
with the Virasoro subalgebra central charge
\bea
c=(n-1)\Biggl{(}1-\frac{n(n+1)}{r(r-1)}\Biggr{)}\ .\nonumber
\ena
The bosonic fields 
$W^{(j)}(z)$ defined via the transform  \eqref{DWAmiura} 
constitute an associative algebra \cite{qWN,FeFr95}. 
For instance, 
the commutation relations between the currents 
$W^{(j)}(z)$ and $W^{(1)}(z)$ 
\bea
W^{(1)}(z)=\sum_{s=1}^{n}\Lambda_s(z)
\label{DWAcur}\ena
are described by 
\bea
&&f^{(j)}(\frac{w}{z})W^{(1)}(z)W^{(j)}(w)-
W^{(j)}(w)W^{(1)}(z)f^{(j)}(\frac{z}{w})=-(x-x^{-1})^2\times
\cr
&&\times
[r]_x[r-1]_x
\{W^{(j+1)}(xw)\delta(x^{j+1}\frac{w}{z})-
W^{(j+1)}(x^{-1}w)\delta(x^{-(j+1)}\frac{w}{z})
\}.\nonumber\\
\label{DWAcom}
\ena
Here $\delta(z)=\sum_{j\in \Z}z^j$ and
\bea
f^{(j)}(z)=\exp{\Biggl(-(x-x^{-1})^2\sum_{m>0}\frac{z^m}{m}
[rm]_x[(r-1)m]_x\frac{[(n-j)m]_x}{[nm]_x}
}
\Biggr). \nonumber
\ena

It is important for us that the relation \eqref{DWAcom} leads 
\cite{FeFr95} to 
\begin{lem} 
\label{LEMDWA1}
The Fourier
modes ${\cal{W}}^{(1)}_t,\ t\in \Z$ of the currents 
$ W^{(1)}(z)$
\bea
{\cal{W}}^{(1)}_t=\oint\frac{dz}{2\pi iz} z^t W^{(1)}(z)
\ena
generate the whole DWA.
\end{lem} 

\subsection{Intertwining property of $\bX_\alpha({\lambda})$.}
The main statement in this section is that 
the screening operators defined by \eqref{SET1}
satisfy the basic property 
of the intertwining operators for DWA:
\begin{prop}
\label{PRDWA1}
\bea
[W^{(j)}(z),\bX_{\alpha}(\lambda)]=0 \ .
\label{DWA2}\ena
\end{prop} 

According to the Lemma \ref{LEMDWA1}  
the proof of the Proposition \ref{PRDWA1} follows
from 
\begin{lem}\label{LEMDWA2}
\bea
[W^{(1)}(z),\bX_{\alpha}(\lambda)]=0\ . 
\label{DWAcom2}
\ena
\end{lem}
We will prove this fact in Appendix \ref{app:3}.
Though we have no proof, we believe that 

\medskip
\noindent
{\bf Conjecture} {\it The screening operators $\eqref{SET1}$
exhaust all the intertwiners for DWA.}
\medskip

\rm 
We expect further that the 
irreducible representations of DWA arise
as the cohomologies of the BRST complex \eqref{BRST}.
\rm

An important remark is that the definition of the DWA 
is invariant under the
change $r\rightarrow 1-r$ \cite{qWN,FeFr95}. 
For this reason there exists another set of 
intertwiners given by the second type
screening operators. 
The construction for them 
is fairly obvious and we will not 
consider this case any more.

Finally, let us stress that the property 
$iv)$ of the Theorem 4.3 which is satisfied 
for $\bX_{\alpha}(\lambda)$ turns out to be
essential in the proof of Lemma \ref{LEMDWA2}.
For this reason arbitrary products of basic 
operators or basic operators itself do not commute 
with DWA generators and can not be treated
as screening operators. 
It has been noted at the end of sec. 2.4 and 6 that
a similar situation takes place in the CFT. 
The general conditions for the existence of 
intertwining operators are discussed in Appendix A.

\setcounter{section}{5}
\setcounter{equation}{0}
\section{CFT limit}
The problem of finding the intertwining operators of 
$W$-algebras \cite{LukFat} 
has been studied for $sl_3$ case (i.e., $n=3$) in \cite{BMP90}, \cite{BMPrev}.
However, it was not clear how to 
generalize the result for arbitrary $n$.
In this section we take the CFT limit $x\rightarrow 1\ , \ i\log(z)\sim 1$
of the basic operators for DWA, and
introduce the notion of basic operators
for W-algebra associated to $sl_n$ algebra with $n\ge2$.
In the case $n=3$ we recover the results of \cite{BMP90}. 

We restrict our discussion to formal level in the sense that
the well-definedness of integrals in the operators is not
considered. In the deformed case, there is no such problem.
The integrals are well-defined because the integrands are single-valued and
the contours are on the unit circle.

Our main goal is to find the formal limit 
of the operators \eqref{SO}
corresponding to the positive root 
$\alpha=\alpha_i+\cdots+\alpha_{i+m}$. Recall that
these basic operators 
have the form:
\bea\label{SOlim}
&&X_\alpha(k_1,\ldots,k_m)=
\oint\cdots\oint\prod_{j=i}^{i+m}\iv{z_j}
{U_i(z_i)}\cdots{U_{i+m}(z_{i+m})}
\nonumber\\
&&\times F(u_i,\ldots,u_{i+m}),
\ena
where the function $F(u_i,\ldots,u_{i+m})$ is
given by
\bea\label{THETlim}
&&F(u_i,\ldots,u_{i+m})=(-1)^{k_1+\cdots+k_m}
\prod_{j=1}^m
{[u_{i+j-1}-u_{i+j}-k_j-\half]\over[u_{i+j-1}-u_{i+j}+\half]}
\nonumber\\
&&
\times\prod_{j=0}^{m}{[u_{i+j}-k_j+k_{j+1}-\half-\hp_{i+j}]
\over [u_{i+j}-\half]}
,
\label{SOS2lim}
\ena
and again $k_0=-1,k_{m+1}=0$ is implied. Let us rewrite 
this expression into the sum of "elementary"
integrals with the specific ordering of the variables on the
unit circle:
\bea
&&
{\displaystyle X_\alpha(k_1,\ldots,k_m)=\sum_{\sigma \in S_{m+1}}\atop }
{{\displaystyle\int\cdots\int}
\atop {\scriptstyle 0<\arg z_i<\cdots <\arg z_{i+m}<2\pi}
}
{\displaystyle\prod_{j=i}^{i+m}\iv{z_j}\atop}
\nonumber\\
&&\times
U_{\sigma(i)}(z_i)\cdots 
{U_{\sigma(i+m)}(z_{i+m})}
F_{\sigma}(u_i,\ldots,u_{i+m})\ ,
\label{SO4lim}
\ena
where 
$S_{m+1}$ is the permutation group of numbers $(i,\ldots, i+m)$
and the function $F_\sigma$ is determined by the condition
\bea
&&U_i(z_{\sigma^{-1}(i)})\cdots U_{i+m}(z_{\sigma^{-1}(i+m)})
F(z_{\sigma^{-1}(i)},\ldots,z_{\sigma^{-1}(i+m)})\\
&=&
U_{\sigma(i)}(z_i)\cdots U_{\sigma(i+m)}(z_{i+m})
F_\sigma(z_i,\ldots,z_{i+m}),\\
\ena
or equivalently by
\bea
F_{\sigma}(u_i,\ldots,u_{i+m})=F(u_{\sigma^{-1}(i)},
\ldots,u_{\sigma^{-1}(i+m)})\prod_{j'< j\atop \sigma(j')-\sigma(j)=1}
(-)\frac{[u_{j}-u_{j'}+\half]}
{[u_{j}-u_{j'}-\half]}\ .
\ena
Now we are able to work out the conformal limit of 
the operators $X_\alpha(k_1,\ldots,k_m)$. 

In what follows in this section, we will use 
the same notations
for the CFT limits of the bosons, screening operators
etc., and mention only
the changes in the definitions. 

First let us discuss the 
limit of the commutation relations for bosons
and the screening current. The first one can be found
directly by setting $x\rightarrow 1$ 
in the formulae \eqref{eqn:2.1}-\eqref{eqn:2.3}.
Namely, 
the oscillators $\beta^j_m$ 
($1\le j\le n-1, m\in\Z\backslash\{0\}$) are defined
by the commutation relations 
\begin{eqnarray}
[\beta^j_m,\beta^k_{m'}]
&=&
m\frac{(n-1)}{n}\frac{(r-1)}{r}\delta_{m+m',0}
\qquad (j=k),
\label{eqn:2.1lim}\\
&=&
-\frac{1}{n}\frac{(r-1)}{r}m\delta_{m+m',0}
\qquad (j\neq k),
\nonumber\\
&&\label{eqn:2.2lim}
\end{eqnarray}
while $\beta^n_m$ is determined via the equation
\begin{equation}
\sum_{j=1}^n\beta^j_m=0.
\label{eqn:2.3lim}
\end{equation}
Note that the definitions of 
the bosonic Fock spaces and 
the zero mode operators $P_\lambda,Q_\lambda$ 
remain the same as in section 2.2. 

The prescription for taking the limit of the 
operator $U_j(z)$ ($j=1,\ldots,n-1$) is also very simple. 
We demand that together with $x\rightarrow 1$, the
parameter $u$ tends to a limiting value in such a way that
$z=x^{2u}$ is fixed. Therefore the screening 
currents of \cite{qWN,FeFr95} 
become the screening currents of \cite{LukFat}:
\begin{eqnarray}
U_j(z)&=&e^{i\sqrt{\frac{r-1}{r}}Q_{\alpha_j}}
z^{\frac{1}{r}\hat{\pi}_{\alpha_j}+\frac{r-1}{r}}
:e^{\sum_{m\neq 0}\frac{1}{m}(\beta^j_m-\beta^{j+1}_m)z^{-m}}:,
\label{eqn:2.4lim}
\end{eqnarray}
where  $z\in \C$.

Note that the product of operators $U(z)U(\zeta)$ is
defined apriori when $|z|>\hskip-3pt>|\zeta|$, and
to be understood as an analytic continuation. When
we compare $U(z)U(\zeta)$ with $U(\zeta)U(z)$ in the deformed case,
there is a common domain of convergence (a neighborhood of the
unit circle $|z/\zeta|=1$). However, in the CFT limit, this will
shrink because poles accumulate to $z/\zeta=1$. Therefore, in order to
compare $U(z)U(\zeta)$ with $U(\zeta)U(z)$ in the CFT limit,
we must specify the path of analytic continuation. In fact, we have
\be
U_i(z)U_j(\zeta)=
\cases{q^{-(\alpha_i,\alpha_j)}U_j(\zeta)U_i(z),&
if $\hbox{\rm arg}\,z<\hbox{\rm arg}\,\zeta$\cr
q^{(\alpha_i,\alpha_j)}U_j(\zeta)U_i(z),&
if $\hbox{\rm arg}\,z>\hbox{\rm arg}\,\zeta$\cr}.
\label{eqn:co1lim}
\en
Here the complex number $q$ is 
\be
q=e^{i\pi {\frac{r-1}{r}}}\ ,
\en
and the left hand side means the analytic continuation from
the region, $\hbox{\rm arg}\,z=\hbox{\rm arg}\,\zeta$ and $|z|>|\zeta|$,
while the right hand side means the analytic continuation from
the region, $\hbox{\rm arg}\,z=\hbox{\rm arg}\,\zeta$ and $|z|<|\zeta|$.

To work out the conformal limit of the theta depending
part of \eqref{SOlim} we parametrize
the variables $x,u$ as following
\bea
&&x=e^{-\epsilon} \quad (\epsilon >0),\cr
&&u=\frac{v}{2i\epsilon} \quad (0<\hbox{\rm Re}\,v<2\pi).
\ena
Now the limit is given by $\epsilon \rightarrow 0$, while
$z=e^{iv}$ remains to be fixed.
According to such a prescription
the function \eqref{THETlim} is changed to be 
the $\PP_{\alpha}$-depending
factor 
\be
\label{THET1lim}
F(u_i,\ldots,u_{i+m})\rightarrow (-1)^{m}q^{-m-k_1-\cdots-k_m}
e^{\frac{i\pi}{r}(1-\hp_i-\cdots-\hp_{i+m})}
\en
if the condition 
$arg(z_{j}) < arg(z_{j+1})$
holds for any $j=i,\ldots,i+m-1$\ .
For a non-trivial transposition $\sigma \in S_{m+1}$, i.e.,
when some of the screening currents $U_{j+1}$
stands to the left of $U_{j}$ the limit has the form
\be
\label{THET2lim}
F_{\sigma}(u_i,\ldots,u_{i+m})\rightarrow q^{f(\sigma)}(-1)^{m}q^{-m-k_1-\cdots-k_m}
e^{\frac{i\pi}{r}(1-\hp_i-\cdots-\hp_{i+m})}
\en
where the function $f(\sigma)$ is 
\be
f(\sigma)=\sum_{j'< j\atop \sigma(j')-\sigma(j)=1}(2k_{\sigma(j)}+1)\ .
\en   
Now it is convenient to introduce the definition 
of an "elementary" integral $I_{i_1\cdots i_{m+1}}$ 
$(1\leq i_1,\ldots,i_{m+1}\leq n-1)$ with the ordering of the 
variables as follows
:
\bea
&&
{\displaystyle I_{i_1\cdots i_{m+1}}=\atop}
{{\displaystyle\int\cdots \int} \atop 
{\scriptstyle 0<arg z_{i}<\cdots <arg z_{i+m}<2\pi}}
{\displaystyle\prod_{j=i}^{i+m}\iv{z_{j}}\atop}
\nonumber\\
&&\times{U_{i_1}(z_i)}\cdots {U_{i_{m+1}}(z_{i+m})}
\ .
\label{SO5lim}
\ena
In terms of these objects the formal limit of the basic
operators \eqref{SO} is given by the expression
\be
X_\alpha(k_1,\ldots,k_m)=q^{-m-k_1-\cdots-k_m}\sum_{\sigma \in S_{m+1}}
q^{f(\sigma)}I_{\sigma(i)\cdots \sigma(i+m)}\ .
\label{BASIClim}
\en
Here, for notational convenience 
we have omitted the irrelevant common
factor
$$
(-1)^me^{\frac{i\pi}{r}(1-\hp_i-\cdots-\hp_{i+m})}\ .
$$
Indeed, whereas the operator $\PP$ itself 
does not commute with screening operators, the important 
properties of \eqref{BASIClim} such as commutations 
with other basic operators and $W$-algebra generators
do not depend on this factor. 

Note
that
the quasi-periodicity 
\eqref{eqn:c}
of basic operators follows from 
\[
q^{r}=(-1)^{r-1}.
\]

The product of two integrals
of the form $I_{i_1\cdots i_m}, I_{i_{m+1}\cdots i_{m+s}}$
is given according to the definition \eqref{SO5lim}
as follows. Let $S_{m+s}$ be the permutation
group of numbers $\{1,\ldots,m+s\}$. Denote by
$S_{m,s}$ the set of elements $\sigma \in S_{m+s}$
such that $\sigma^{-1}(j)<\sigma^{-1}(j+1)$ for each
$j\neq m,m+s$.
Then 
\bea
I_{i_1\cdots i_m}I_{i_{m+1}\cdots i_{m+s}}=
\sum_{\sigma \in S_{m,s}}
q^{g_{m,s}(\sigma)}I_{i_{\sigma(1)}\cdots i_{\sigma(m+s)}
}
\label{PRODcft}
\ena
where
\be
g_{m,s}(\sigma)=\sum_{j<l \atop \sigma(j)>\sigma(l)}
(\alpha_{i_{\sigma(j)}},\alpha_{i_{\sigma(l)}})
\en
For instance, 
let us 
work out the decomposition of 
the product of two basic
operators (see also \cite{Fel89,BMP90}):
\bea
&&X_{12}(k)=q^{-k-1}I_{12}+q^{k}I_{21} ,\\
&&X_{1}=I_1,
\ena
into the "elementary" integrals.
In this example we will also show the
validity of \eqref{eqn:a3}. 
The formal product of these operators
is the operator of weight $-2\alpha_1-\alpha_2$. It
can be derived from the definition \eqref{PRODcft}
\bea
&&I_{i_1i_2}I_{i_3}=I_{i_1i_2i_3}
+q^{(\alpha_{i_3},\alpha_{i_2})}I_{i_1i_3i_2}+
q^{(\alpha_{i_3},\alpha_{i_1})+(\alpha_{i_3},\alpha_{i_2})}I_{i_3i_1i_2},
\\
&&
I_{i_1}I_{i_2i_3}=I_{i_1i_2i_3}
+q^{(\alpha_{i_1},\alpha_{i_2})}I_{i_2i_1i_3}+
q^{(\alpha_{i_1},\alpha_{i_3})+(\alpha_{i_1},\alpha_{i_2})}I_{i_2i_3i_1}.
\ena
In particular,
\bea
&&I_{12}I_1=(q+q^{-1})I_{112}+I_{121},
\nonumber\\
&&I_{21}I_1=qI_{121}+(1+q^2)I_{211},
\nonumber\\
&&I_1I_{12}=(1+q^2)I_{112}+qI_{121},
\nonumber\\
&&I_1I_{21}=I_{121}+(q+q^{-1})I_{211}.
\ena
Using these formulae we find that
\bea
&&X_1X_{12}(k)
\nonumber\\
&&=q^{-k}(q+q^{-1})I_{112}+(q^{k}+q^{-k})I_{121}+
+q^{k}(q+q^{-1})I_{211}
\nonumber\\
&&=X_{12}(k-1)X_1.
\label{SL3lim}
\ena
This confirms that 
\eqref{eqn:a3} still holds in the CFT limit.
Similarly, one can 
check \eqref{eqn:a1},\eqref{eqn:a2},\eqref{eqn:a4} 
and similar
properties of the basic operators in the $n>3$ case.

It would be useful, however, to express the basic operators
in the "conventional" form, i.e., as non-commutative 
polynomials of  
$X_j$. Let us first prepare the necessary notations.
Introduce the
bracket
\bea
&&\{A,B\}_k\equiv-[k]_qAB+[k+1]_qBA\ ,
\label{K1lim}\\
&&\{A,B\}_0=BA\ ,
\label{K2lim}\\
&&\{A,B\}_{-1}=AB\ ,
\label{K3lim}\\
&&\{\{A,B\}_{k_1},C\}_{k_2}=\{A,\{B,C\}_{k_2}\}_{k_1}\ , if [A,C]=0 .
\label{K4lim}
\ena
Now one finds that the following Lemma 
holds:
\begin{lem} 
\label{LIM}
\bea
&&q^{-m-k_1-\cdots-k_m}\sum_{\sigma \in S_{m+1}}
q^{f(\sigma)}I_{\sigma(i)\cdots \sigma(i+m)}
\nonumber\\
&&=q^{1-m-k_1-\cdots-k_{m-1}}\sum_{\sigma \in S_{m}}
q^{f(\sigma)}\{I_{\sigma(i)\cdots \sigma(i+m-1)},X_m\}_{k_{m}}\ .
\label{TRANSlim}
\ena
\end{lem} 
The proof follows from a straightforward 
decomposition of the right hand side into 
"elementary integrals" $I_{i_1\cdots i_{m+1} }$.
Applying the equation \eqref{TRANSlim}
we arrive at the 
\begin{prop}
\label{LIM2}
\be
\label{SOcft}
X_{i\cdots i+m}(k_1,\ldots,k_m)=
\{\ldots\{\{X_i,X_{i+1}\}_{k_1},X_{i+2}\}_{k_2},\ldots,X_{i+m}\}_{k_m}\ . 
\en
\end{prop}
In particular, a conformal analogue of the operators \eqref{eqn:X12}
is
\begin{equation}
X_{12}(k)=-[k]_qX_1X_2+[k+1]_qX_2X_1.
\label{cft1}
\end{equation}
Note that using this representation the equality of
two decompositions \eqref{SL3lim} and CFT analogue 
for \eqref{eqn:a4}
can be rewritten
as $q$-Serre relations \cite{BMP90}    
\begin{eqnarray*}
&&[k]_qX_{j}^2X_{j+1}-([k+1]_q+[k-1]_q)
X_{j}X_{j+1}X_{j}+[k]_qX_{j+1}X_{j}^2=0 ,
\\
&&[k]_qX_{j+1}^2X_{j}-([k+1]_q+[k-1]_q)
X_{j+1}X_{j}X_{j+1}+[k]_qX_{j}X_{j+1}^2=0 .
\end{eqnarray*}
These equations together with the 
commutativity of screening currents
$
U_{j},\ U_{l}
$ for $|j-l|\neq 1$
imply that
in the CFT limit, the screening operators corresponding to the
simple roots satisfy the relations for the nilpotent half of 
the quantum group $U_q(sl_n)$. 

Using the properties of the bracket \eqref{K1lim}-\eqref{K4lim} 
one can check that 
the operators \eqref{SOcft} satisfy the basic relations \eqref{QR1}-\eqref{QR7}
of Lemma B.4. which now have the form:
\bea
&&X_{i\cdots i+l-1}(k_1,\ldots,k_{l-1})X_{i+l\cdots i+l+m}
(k_{l+1},\ldots,
k_{l+m})
\nonumber\\
&&
=X_{i\cdots i+l+m}(k_1,\ldots,
k_{l-1},-1,k_{l+1},\ldots,k_{l+m}),
\nonumber\\\label{cftQR1}\\
&&X_{i+l\cdots i+l+m}(k_{l+1},\ldots,k_{l+m})X_{i\cdots i+l-1}(k_1,\ldots,k_{l-1})
\nonumber\\
&&
=X_{i\cdots i+l+m}(k_1,\ldots,k_{l-1},
0,k_{l+1},\ldots,k_{l+m}),
\nonumber\\\label{cftQR2}\\
&&X_iX_{i\cdots i+m}(k_1,k_2,\ldots,k_m)
=X_{i\cdots i+m}(k_1-1,k_2,\ldots,k_m)X_i,
\nonumber\\
\label{cftQR21}\\
&&X_{i+m}X_{i\cdots i+m}(k_1,\ldots,k_{m-1},k_m-1)=
X_{i\cdots i+m}(k_1,\ldots,k_{m-1},k_m)X_{i+m},
\nonumber\\
\label{cftQR6}\\
&&X_{i\cdots i+m}(k_1,\ldots,k_m)X_{j\cdots j+l}(k'_1,\ldots,k'_l)=
X_{j\cdots j+l}(k'_1,\ldots,k'_l)
\nonumber\\
&&
\times X_{i\cdots i+m}(k_1,\ldots,k_m),
\phantom{mmmmmmm}\phantom{mmmmmmm}\hbox{ if $i+m+1<j$.}
\nonumber\\
\label{cftQR7}
\ena

In the case $n=3$ the procedure given 
in section 3 may be followed in CFT limit to
obtain the intertwining operators for
$W_3$ algebra starting from the basic
operators \eqref{cft1}. It can be easily verified that
the expression for $\bX_{\alpha}(\lambda)$ constructed
in such a way coincides with the result of \cite{BMP90}. 
We also inspect directly in CFT limit
the commutativity property of basic 
operators for the
$n=4$.

We believe that for general $n$, the construction of the intertwining
operators $\bX_{\alpha}(\lambda)$ for
$W$ algebra is identical for those for the DWA \eqref{MSC} with
the only difference in the
definition of the basic operators \eqref{SOcft}.

Our discussion was rather formal since 
we did not examine the well-defineness of the
operators. 
To treat $\bX_{\alpha}(\lambda)$ as operators on the Fock space one
must fix the contour prescription \cite{BMP90}.
We do not discuss this problem in this paper.

\setcounter{section}{6}
\setcounter{equation}{0}
\section{Discussions}
In this paper, we constructed the intertwining operators
which commute with the action of the deformed $W_n$-algebra
on the bosonic Fock spaces \cite{qWN,FeFr95}. In the conformal limit, the
case we have discussed in this paper corresponds to 
representations of the $W_n$-algebra with the central charge
\bea
c&=&(n-1)\Biggl(1-{n(n+1)\over r(r-1)}\Biggr)
\ena
where $r$ is an integer such that $r\ge n+2$.
In the language of solvable lattice models,
they correspond to the $sl_n$ RSOS models \cite{RSOS}.

The main difference in our construction compared to the
corresponding conformal limit, is that we need to construct basic operators
which change the weight of Fock spaces by a positive root.
In the conformal case, the basic operators can be expressed in terms
of the operators corresponding to the simple roots. In other words,
the case $n\ge3$ contains a new feature which was not seen
in the $n=2$ case considered in \cite{LukPug2}.

One can look at this situation in the following way.
It is well-known that the algebra of the screening operators
in the conformal field theory is isomorphic to the
algebra $U_q(b_+)$ with $q=e^{\pi i(r-1)/r}$.
An elliptic deformation of this algebra was considered in \cite{FeOd95}.
In this paper, we identified it with the algebra of the
screening operators of the deformed $W_n$-algebra, and we derived
a set of quadratic relations among the generators of that algebra.
These relations can be considered as the elliptic deformation
of the Serre relations.

One of our original aims was to construct the Felder-type
complex for the irreducible $W_n$ modules. This was not achieved in
two reasons. First of all, except for the case of $n=3$,
we cannot fix the signs so that we have $d^2=0$ for the coboundary operators.
Secondly, even for the case $n=3$, we have no result on the cohomology of the
complex.

We considered the CFT limit of our construction. Our discussion stays
in a formal level because we only considered the limit of
the integrands of the screening operators.

\bigskip
{\bf Acknowledgment.} \quad 
B.F. and A. O. would like to thank 
RIMS for the hospitality throughout their stay, during which 
this work was done.  T.M. thanks IHP and ENS where
he stayed in the last stage of the work.
This work is partially 
supported by Grant-in-Aid for Scientific Research on Priority 
Areas 231, the Ministry of Education, Science and Culture.

\appendix
\setcounter{equation}{0}
\section{Admissibility and commuting squares}\label{app:1}

\subsection{Admissible pairs}

Here we derive the condition for the admissibility of a pair 
$(\lambda,\lambda^\alpha)$. We follow the notation of Section \ref{SEC3}.
In particular, we suppose $\lambda=t_\gamma\sigma\Lambda$ throughout this
section.

\begin{lem}\label{LEM1}
The condition $(\sigma\Lambda,\alpha)\gl0$ is equivalent to
$l(r_\alpha\sigma)\gl l(\sigma)$.
\end{lem}
\proof
Suppose that $r_\alpha=(ij)$ is the transposition and
\bea
\sigma=\bigl(\sigma(1),\ldots,i,\ldots,j,\ldots,\sigma(n)\bigr).
\ena
We have
\bea
r_\alpha\sigma=\bigl(\sigma(1),\ldots,j,\ldots,i,\ldots,\sigma(n)\bigr).
\ena
Therefore 
$l(r_\alpha\sigma)\gl l(\sigma)$
is equivalent to $i\lg j$.
Since $\alpha$ is positive and $r_\alpha=(ij)$, we have
$\alpha=\ve_i-\ve_j$ if $i<j$ and $\alpha=\ve_j-\ve_i$ if $i>j$.
Noting that $\sigma^{-1}\ve_i=\ve_{\sigma^{-1}(i)}$
and $\sigma^{-1}(i)<\sigma^{-1}(j)$,
we conclude that
$l(r_\alpha\sigma)\gl l(\sigma)$ is equivalent to
$\sigma^{-1}\alpha\gl0$, and therefore to
$(\sigma\Lambda,\alpha)=(\Lambda,\sigma^{-1}\alpha)\gl0$.
\qed

\begin{lem}
$l(r_\alpha)=2|\alpha|-1$.
\end{lem}
\proof
If $\alpha=\alpha_{i\cdots i+m}$, a reduced expression of
$r_\alpha$ is given by
\bea
r_\alpha=s_i\cdots s_{i+m-1}s_{i+m}s_{i+m-1}\cdots s_i.\nonumber
\ena
\qed

We consider an operator $X_\alpha(\lambda)$ if and only if
\bea
d_\alpha(\lambda)
\buildrel\rm def\over=
\deg(\lambda^\alpha)-\deg(\lambda)=1.
\ena
Note that
\bea
d_\alpha(\lambda)
=\cases{l(r_\alpha\sigma)-l(\sigma)& if $(\sigma\Lambda,\alpha)>0$;\cr
l(r_\alpha\sigma)-l(\sigma)+2|\alpha|& if $(\sigma\Lambda,\alpha)<0$.\cr}
\ena
In particular, we have
\bea\label{INP}
0<d_\alpha(\lambda)<2|\alpha|.
\ena

\begin{lem}\label{DL}
A pair $(\lambda,\lambda^\alpha)$ is admissible if and only if
one of the following is holds:
\par
{\rm (i)} $(\sigma\Lambda,\alpha)>0$, and
$(\sigma\Lambda,\beta)<0$ or 
$(\sigma\Lambda,\gamma)<0$ for any partition $\alpha=\beta+\gamma$
$(\beta,\gamma\in\Delta_+)$.
\par
{\rm (ii)} $(\sigma\Lambda,\alpha)<0$, and
$(\sigma\Lambda,\beta)<0$ and 
$(\sigma\Lambda,\gamma)<0$ for any partition $\alpha=\beta+\gamma$
$(\beta,\gamma\in\Delta_+)$.
\end{lem}
In particular, $(\lambda,\lambda^\alpha)$ is always admissible 
for a simple root $\alpha=\alpha_j$. 

\proof
We follow the argument in the proof of Lemma \ref{LEM1}.
If $l(r_\alpha\sigma)>l(\sigma)$, we set
$\beta=\ve_i-\ve_k$ and $\gamma=\ve_k-\ve_j$ for $k$ such that $i<k<j$.
The condition $d_\alpha(\lambda)=l(r_\alpha\sigma)-l(\sigma)=1$
is equivalent to $\sigma^{-1}(k)<\sigma^{-1}(i)$
or $\sigma^{-1}(j)<\sigma^{-1}(k)$ for any such $k$.
This is equivalent to $(\sigma\Lambda,\beta)<0$ or
$(\sigma\Lambda,\gamma)<0$, respectively.
If $l(r_\alpha\sigma)<l(\sigma)$, we set
$\beta=\ve_j-\ve_k$ and $\gamma=\ve_k-\ve_i$ for $k$ such that $j<k<i$.
The condition $d_\alpha(\lambda)=l(r_\alpha\sigma)-l(\sigma)+2|\alpha|=1$
is equivalent to $\sigma^{-1}(i)<\sigma^{-1}(k)<\sigma^{-1}(j)$
for any such $k$. This is equivalent to $(\sigma\Lambda,\beta)<0$ and
$(\sigma\Lambda,\gamma)<0$.
\qed

\subsection{Commuting squares}

In this subsection we prove Theorem \ref{thm:comsq}. 

The assertion (i) follows immediately from 
\bea
d_\alpha(\lambda)+d_\alpha(\lambda^\alpha)=2|\alpha|.
\ena
Below we shall prove the assertion (ii) case-by-case. 

\medskip

\noindent{Case $(\alpha,\beta)=0$}

Set $m=m_\alpha(\lambda)$ and $m'=m_\beta(\lambda)$. 
Recall (\ref{LA}). 
We have
\bea
(r_\alpha\sigma\Lambda,\beta)=(\sigma\Lambda,\beta).
\ena
This implies $m_\beta(\lambda^\alpha)=m'$. Therefore, we have
\bea
\lambda-\lambda^{\alpha,\beta}=m\alpha+m'\beta.
\label{AAA}
\ena
Similarly, we have
\bea
\lambda-\lambda^{\beta,\alpha}=m\alpha+m'\beta.
\ena
This implies
\bea
d_\beta(\lambda)+d_\alpha(\lambda^\beta)&=&
\deg(\lambda^{\beta,\alpha})-\deg(\lambda)\nonumber\\
&=&\deg(\lambda^{\alpha,\beta})-\deg(\lambda)\nonumber\\
&=&2.
\ena
Since $d_\beta(\lambda),d_\alpha(\lambda^\beta)>0$, we have
\bea
d_\beta(\lambda)=d_\alpha(\lambda^\beta)=1.\label{BBB}
\ena
Namely, $(\lambda,\lambda^\beta,\lambda^{\beta,\alpha})$
is admissible.

Let us show the uniqueness of $\alpha',\beta'$. 
Suppose that $\alpha=\ve_i-\ve_j$ and $\beta=\ve_k-\ve_l$.
We consider only the case when $k<i<j<l$ and set $\gamma_1=\ve_k-\ve_i$
and $\gamma_2=\ve_j-\ve_l$. The other cases are similar.
If
\bea
(\lambda,\lambda^\gamma,\lambda-m\alpha-m'\beta)
\ena
is admissible and $\gamma\not=\alpha,\beta$, then we have
\bea
\gamma=\ve_k-\ve_j=\gamma_1+\alpha\hbox{ \rm or }
\gamma=\ve_i-\ve_l=\alpha+\gamma_2
\ena
and
\bea
m=m'=m_{\gamma_1+\alpha}(\lambda)=m_{\alpha+\gamma_2}(\lambda).
\ena
Since $m_\alpha(\lambda)\equiv (\sigma \Lambda,\alpha)$ and 
$m_{\gamma_1+\alpha}(\lambda)\equiv (\sigma \Lambda,\gamma_1+\alpha)$ 
$\bmod r$, we have 
$(\sigma\Lambda,\gamma_1)\equiv 0\bmod r$. 
This is a contradiction. 

\medskip

\noindent{Case $(\alpha,\beta)=1$}

Set $m=m_\alpha(\lambda)$ and $m'=m_\beta(\lambda^\alpha)$.
We have (\ref{AAA}). The only way other than (\ref{AAA})
to write $\lambda-\lambda^{\alpha,\beta}$ as a positive linear combination
of two positive roots is
\bea\label{ONLY}
\lambda-\lambda^{\alpha,\beta}
=\cases{m(\alpha-\beta)+(m+m')\beta& if $\alpha-\beta\in\Delta_+$;\cr
(m+m')\alpha+m'(\beta-\alpha)& if $\beta-\alpha\in\Delta_+$.\cr}
\ena
The uniqueness is then obvious from (\ref{ONLY}). 
Note that
\bea
&&(r_\beta\sigma\Lambda,\alpha-\beta)
=(\sigma\Lambda,\alpha)=\cases{m& if $(\sigma\Lambda,\alpha)>0$;\cr
m-r& if $(\sigma\Lambda,\alpha)<0$,\cr}\label{E1}\\
&&(r_\alpha\sigma\Lambda,\beta)
=(\sigma\Lambda,\beta-\alpha)=\cases{
m'& if $(\sigma\Lambda,\alpha-\beta)<0$;\cr
m'-r& if $(\sigma\Lambda,\alpha-\beta)>0$.\cr}
\label{E2}
\ena

Let us show that $(\lambda,\lambda^\beta,\lambda^{\beta,\alpha-\beta})$
is admissible if $\alpha-\beta\in\Delta_+$.
From (\ref{E1}) follows $m_{\alpha-\beta}(\lambda^\beta)=m$.
Let us prove
\bea
m_\beta(\lambda)=m+m'.\label{CC}
\ena
If $(\sigma\Lambda,\alpha)>0$
and $(\sigma\Lambda,\alpha-\beta)<0$, the statement (\ref{CC})
follows from (\ref{E1}) and (\ref{E2}).
The case $(\sigma\Lambda,\alpha)<0$
and $(\sigma\Lambda,\alpha-\beta)>0$ contradicts with Lemma \ref{DL}
because $(\lambda,\lambda^\alpha)$ is admissible.
In the remaining cases, we have
$(\sigma\Lambda,\beta)=m+m'-r$. From Lemma \ref{DL} (applied to
$(\lambda,\lambda^\alpha)$), we have $m+m'-r<0$,
and therefore (\ref{CC}).
Now, we have $\lambda^{\beta,\alpha-\beta}=\lambda^{\alpha,\beta}$.
This implies (see (\ref{BBB}))
\bea
d_\alpha(\lambda)=d_{\alpha-\beta}(\lambda^\beta)=1.
\ena
Thus we proved the admissibility of
$(\lambda,\lambda^\beta,\lambda^{\beta,\alpha-\beta})$
if $\alpha-\beta\in \Delta_+$.

Next, we show that
$(\lambda,\lambda^{\beta-\alpha},\lambda^{\beta-\alpha,\alpha})$
is admissible if $\beta-\alpha\in\Delta_+$.
From (\ref{E2}) follows $m_{\beta-\alpha}(\lambda)=m'$.
Let us prove
\bea
m_\alpha(\lambda^{\beta-\alpha})=m+m'.\label{BB}
\ena
Note that 
$(r_{\beta-\alpha}\sigma\Lambda,\alpha)=(\sigma\Lambda,\beta)$.
If $(\sigma\Lambda,\alpha)>0$ and $(\sigma\Lambda,\beta-\alpha)>0$,
the statement (\ref{BB}) follows
from (\ref{E1}) and (\ref{E2}). 
The case $(\sigma\Lambda,\alpha)<0$ and
$(r_\alpha\sigma\Lambda,\beta)=(\sigma\Lambda,\beta-\alpha)<0$
contradicts with Lemma \ref{DL} because
$(\lambda^\alpha,\lambda^{\alpha,\beta})$ is admissible
and $(r_\alpha\sigma\Lambda,\alpha)=-(\sigma\Lambda,\alpha)>0$.
In the remaining cases, we have $(\sigma\Lambda,\beta)=m+m'-r$.
From Lemma \ref{DL} (applied to $(\lambda^\alpha,\lambda^{\alpha,\beta})$)
we have $(r_\alpha\sigma\Lambda,\beta-\alpha)=(\sigma\Lambda,\beta)<0$,
and therefore (\ref{BB}).
Thus we proved (ii) when $(\alpha,\beta)=1$.

\medskip

\noindent{Case $(\alpha,\beta)=-1$}

Set $m=m_\alpha(\lambda)$ and $m'=m_\beta(\lambda^\alpha)$. 
Because of (\ref{3.00}) we have $(\sigma\Lambda,\beta)\not\equiv0\bmod r$,
and therefore $m\not=m'$. We have (\ref{AAA}).
The only way other than (\ref{AAA}) to write
$\lambda-\lambda^{\alpha,\beta}$ as a positive linear combination of two
positive roots is
\bea\label{ONLY2}
\lambda-\lambda^{\alpha,\beta}=
\cases{(m-m')\alpha+m'(\alpha+\beta)& if $m>m'$;\cr
m(\alpha+\beta)+(m'-m)\beta& if $m<m'$.\cr}
\ena
Again he uniqueness is obvious from (\ref{ONLY2}). 
Note that
\bea
&&(r_\beta\sigma\Lambda,\alpha+\beta)=
(\sigma\Lambda,\alpha)=\cases{m& if $(\sigma\Lambda,\alpha)>0$;\cr
m-r& if $(\sigma\Lambda,\alpha)<0$,\cr}\label{E5}\\
&&(r_\alpha\sigma\Lambda,\beta)=(\sigma\Lambda,\alpha+\beta)=
\cases{m'& if $(\sigma\Lambda,\alpha+\beta)>0$;\cr
m'-r& if $(\sigma\Lambda,\alpha+\beta)<0$.\cr}\label{E6}
\ena

Let us show that
$(\lambda,\lambda^{\alpha+\beta},\lambda^{\alpha+\beta,\alpha})$
is admissible if $m>m'$.
From (\ref{E6}) follows $m_{\alpha+\beta}(\lambda)=m'$.
Let us prove
\bea
m_\alpha(\lambda^{\alpha+\beta})=m-m'.\label{AB}
\ena
Note that
$(r_{\alpha+\beta}\sigma\Lambda,\alpha)=-(\sigma\Lambda,\beta)$.
If $(\sigma\Lambda,\alpha)\gl0$ and $(\sigma\Lambda,\alpha+\beta)\gl0$,
we have $-(\sigma\Lambda,\beta)=m-m'>0$,
and therefore (\ref{AB}).
If $(\sigma\Lambda,\alpha)<0$ and $(\sigma\Lambda,\alpha+\beta)>0$,
we have $-(\sigma\Lambda,\beta)=m-m'-r<0$,
and therefore (\ref{AB}). If $(\sigma\Lambda,\alpha)>0$
and $(\sigma\Lambda,\alpha+\beta)<0$, we have
$-(\sigma\Lambda,\beta)=m-m'+r>r$. This is a contradiction.

Let us show that
$(\lambda,\lambda^\beta,\lambda^{\beta,\alpha+\beta})$
is admissible if $m<m'$. From (\ref{E6}) we have
$m_{\alpha+\beta}(\lambda^\beta)=m$. Let us prove
\bea
m_\beta(\lambda)=m'-m.\label{BA}
\ena
If $(\sigma\Lambda,\alpha)\gl0$ and $(\sigma\Lambda,\alpha+\beta)\gl0$,
we have $(\sigma\Lambda,\beta)=m'-m>0$,
and therefore (\ref{BA}).
If $(\sigma\Lambda,\alpha)<0$ and $(\sigma\Lambda,\alpha+\beta)>0$,
we have $(\sigma\Lambda,\beta)=m'-m+r>r$. This is a contradiction.
If $(\sigma\Lambda,\alpha)>0$ and $(\sigma\Lambda,\alpha+\beta)<0$,
we have $(\sigma\Lambda,\beta)=m'-m-r<0$,
and therefore (\ref{BA}).
We have completed the proof of (ii) when $(\alpha,\beta)=-1$.

\setcounter{equation}{0}
\section{Generalized Serre relations}\label{app:2}

We modify the relations (\ref{eqn:co1}), (\ref{eqn:co2}) and (\ref{ZMCOM})
(keeping (\ref{eqn:co3})) as follows.
\bea
\xi_i(u)\xi_i(v)&=&{[u-v-\delta]\over[u-v+\delta]}\xi_i(v)\xi_i(u),\label{MR1}\\
\xi_i(u)\xi_j(v)&=&{[u-v+{\delta\over2}]\over[u-v-{\delta\over2}]}\xi_j(v)\xi_i(u)
\hbox{ if $|i-j|=1$},\label{MR2}\\
\hp_i\xi_j(u)&=&\xi_j(u)\bigl(\hp_i-(\alpha_i,\alpha_j)\delta\bigr).\label{MR4}
\ena
Here, $\delta$ is a parameter. Note that if we set $\delta=0$
we have a commutative algebra.

Fix
\bea
(a_1,\ldots,a_{n-1})\quad(a_i\in\Z_{\geq0}).\label{3.0}
\ena

Consider a function $f$ of the
variables $u^{(b)}_j$ $(1\leq j\leq n-1,1\leq b\leq a_j)$ and $\kappa_j$
$(1\leq j\leq n-1)$. We assume that $f$ is symmetric in
$(u^{(1)}_j,\ldots,u^{(a_j)}_j)$ for each $1\leq j\leq n-1$.
We call $f$ a function of type $(a_1,\ldots,a_{n-1})$.

Suppose that $f$ is of type $(a_1,\ldots,a_{n-1})$ and
$g$ is of type $(b_1,\ldots,b_{n-1})$. We define the $*$-product
$f*g$ of $f$ and $g$
to be the function of type $(a_1+b_1,\ldots,a_{n-1}+b_{n-1})$
given by
\bea
&&(f*g)(u^{(1)}_1,\ldots,u^{(a_1)}_1,
v^{(1)}_1,\ldots,v^{(b_1)}_1,\ldots,
u^{(1)}_{n-1},\ldots,u^{(a_{n-1})}_{n-1},
v^{(1)}_{n-1},\ldots,
v^{(b_{n-1})}_{n-1};
\nonumber\\
&&
\kappa_1,\ldots,\kappa_{n-1})
=S\Biggl(
f(u^{(1)}_1,\ldots,u^{(a_1)}_1,\ldots,u^{(1)}_{n-1},\ldots,
u^{(a_{n-1})}_{n-1};
\kappa_1+(-2b_1+b_2)\delta,
\nonumber\\
&&
\kappa_2+(b_1-2b_2+b_3)\delta,
\ldots,\kappa_{n-1}+(b_{n-2}-2b_{n-1})\delta)\nonumber\\
&&\times
g(v^{(1)}_1,\ldots,v^{(b_1)}_1,\ldots,v^{(1)}_{n-1},\ldots,v^{(b_{n-1})}_{n-1};
\kappa_1,\kappa_2,\ldots,\kappa_{n-1})\nonumber\\
&&\prod_{{1\leq j\leq n-1\atop1\leq a\leq a_j}\atop1\leq b\leq b_j}
{[u^{(a)}_j-v^{(b)}_j-\delta]\over[u^{(a)}_j-v^{(b)}_j]}
\prod_{{1\leq j\leq n-2\atop 1\leq a\leq a_j}\atop1\leq b\leq b_{j+1}}
{[u^{(a)}_j-v^{(b)}_{j+1}+{\delta\over2}]\over[u^{(a)}_j-v^{(b)}_{j+1}]}
\prod_{{2\leq j\leq n-1\atop 1\leq a\leq a_j}\atop1\leq b\leq b_{j-1}}
{[u^{(a)}_j-v^{(b)}_{j-1}+{\delta\over2}]\over[u^{(a)}_j-v^{(b)}_{j-1}]}
\Biggr).
\nonumber\\
\ena
Here the symbol $S$ means the symmetrization of
$(u^{(1)}_j,\ldots,u^{(a_j)}_j,v^{(1)}_j,\ldots,v^{(b_j)}_j)$
for each $1\leq j\leq n-1$.

Let $f_{i\cdots i+m}^{(k_1,\ldots,k_m)}$ be the following function of type
$(a_1,\ldots,a_{n-1})$ with
\bea
&&a_j=\cases{1& if $i\leq j\leq i+m$;\cr
0& otherwise,\cr}\nonumber\\
&&f^{(k_1,\ldots,k_m)}_{i\cdots i+m}
(u_i,\ldots,u_{i+m};\kappa_i,\ldots,\kappa_{i+m})
=\prod_{j=1}^m{[u_{i+j-1}-u_{i+j}-(k_j+{1\over2})\delta]\over
[u_{i+j-1}-u_{i+j}]}\nonumber\\
&&\qquad \times\prod_{j=0}^{m}
[u_{i+j}-(k_j-k_{j+1}+{1\over2})\delta-\kappa_{i+j}]
\qquad (k_0=-1,k_{m+1}=0).
\label{SOS}
\ena
If $m=0$, we understand the function $f_i$ of type
$(0,\ldots,0,{\buildrel {i\atop\vee} \over 1},0,\ldots,0)$,
to be $[u^{(1)}_i+{\delta\over2}-\kappa_i]$.

\begin{thm}
\bea
f_{i\cdots i+m}^{(k_1,\ldots,k_m)}*f_{i\cdots i+m}^{(k_1+p,\ldots,k_m+p)}
=f_{i\cdots i+m}^{(k_1+p,\ldots,k_m+p)}*f_{i\cdots i+m}^{(k_1,\ldots,k_m)}.
\ena
\end{thm}
\proof
This is similar to Theorem \ref{COM}. \qed

Set
\bea
f_{i\dots i+m}[k_1,\ldots,k_m]&=&f_{i\cdots i+m}^{(k_1,\ldots,k_m)},\\
f^{(a)}_{i\cdots i+m}[k_1\ldots,k_m]&=&*\prod_{b=1}^a
f_{i\cdots i+m}[k_1-b+1,\ldots,k_m-b+1]
\ena
where the symbol $*$ in front of the usual product symbol means that this is
a $*$-product. The functions $f^{(a)}_{i\cdots i+m}[k_1,\ldots,k_m]$
satisfy a set of quadratic relations in $*$-product. They are given below.
By specialization $\delta=1$, we get the relations for the
screening operators (\ref{MSC}) $X^{(a)}_{i\cdots i+m}(k_1,\ldots,k_m)$.

For the proof of the quadratic relations we prepare a lemma.

Let $F$ be the algebra over $\C$ with the $*$-product, that is generated
by elements $f_{i\cdots i+m}[k_1,\ldots,k_m]$. The algebra $F$
is graded
\bea
F&=&\oplus_{(a_1,\ldots,a_{n-1})\in\Z^{n-1}_{\geq0}}
F_{a_1,\ldots,a_{n-1}}
\ena
where $F_{a_1,\ldots,a_{n-1}}$ consists of the functions of type
$(a_1,\ldots,a_{n-1})$.

\begin{lem}
If $f,g\in F$ and $f*g=0$, then $f=0$ or $g=0$.
\end{lem}

\proof
Suppose that $f$ is of type $(a_1,\ldots,a_{n-1})$ and
$g$ is of type $(b_1,\ldots,b_{n-1})$.
We expand $f$ and $g$ in power series of $\delta$.
If both $f$ and $g$ are non-zero, we have
\bea
f&=&f_0\delta^{m_1}+o(\delta^{m_1}),\quad f_0\not=0,\nonumber\\
g&=&g_0\delta^{m_2}+o(\delta^{m_2}),\quad g_0\not=0,\nonumber
\ena
for some $m_1$ and $m_2$. From $f*g=0$ follows that
\bea
&&S\Bigl(f_0(u^{(1)}_1,\ldots,u^{(a_1)}_1,
\ldots,u^{(1)}_{n-1},\ldots,u^{(a_{n-1})}_{n-1};
\kappa_1,\ldots,\kappa_{n-1})\nonumber\\
&&\cdot g_0(v^{(1)}_1,\ldots,v^{(b_1)}_1,
\ldots,v^{(1)}_{n-1},\ldots,v^{(b_{n-1})}_{n-1};
\kappa_1,\ldots,\kappa_{n-1})\Bigr)=0\nonumber\\
\ena
We will show that $f_0=0$ or $g_0=0$. Choose
$w_1,\ldots,w_{n-1}\in\C$ so that $f_0$ and $g_0$ are holomorphic at
$u^{(b)}_j=w_j$ and $v^{(b)}_j=w_j$, respectively.
Power series expansion in $u^{(b)}_j-w_j$ and $v^{(b)}_j-w_j$
reduces the problem to the case when $f_0$ and $g_0$
are symmetric polynomials in
$(u^{(1)}_j,\ldots,u^{(a_j)}_j)$
and $(v^{(1)}_j,\ldots,v^{(a_j)}_j)$, respectively.
Finally the following lemma reduces the problem to the case of the polynomial
ring.

\begin{lem}
Let $G_{a_1,\ldots,a_{n-1}}$ be the $\C$-linear space of
polynomials in the variables
\[
u^{(1)}_1,\ldots,u^{(a_1)}_1,\ldots,
u^{(1)}_{n-1},\ldots,u^{(a_{n-1})}_{n-1},
\]
that are symmetric in $(u^{(1)}_j,\ldots,u^{(a_j)}_j)$
for each $1\leq j\leq n-1$.
Set
\bea
G=\oplus_{(a_1,\ldots,a_{n-1})\in\Z^{n-1}_{\geq0}}
G_{a_1,\ldots,a_{n-1}}.
\ena
Define the $*$-product in $G$ by
\bea
&&f\in G_{a_1,\ldots,a_{n-1}},
g\in G_{b_1,\ldots,b_{n-1}}\rightarrow
f*g\in G_{a_1+b_1,\ldots,a_{n-1}+b_{n-1}},\nonumber
\ena
where
\bea
&&(f*g)(u^{(1)}_1,\ldots,u^{(a_1)}_1,v^{(1)}_1,\ldots,v^{(b_1)}_1,
\ldots,u^{(1)}_{n-1},\ldots,
u^{(a_{n-1})}_{n-1},v^{(1)}_{n-1},\ldots,v^{(b_{n-1})}_{n-1})\nonumber\\
&&=
S\Bigl(f(u^{(1)}_1,\ldots,u^{(a_1)}_1,
\ldots,u^{(1)}_{n-1},\ldots,u^{(a_{n-1})}_{n-1})
g(v^{(1)}_1,\ldots,v^{(b_1)}_1,
\ldots,v^{(1)}_{n-1},\ldots,
\nonumber\\
&&\ldots,v^{(b_{n-1})}_{n-1})\Bigr).\nonumber
\ena
There is a ring homomorphism between $G$ and the polynomial ring of the
variables
\[
(x^{(0)}_1,x^{(1)}_1,x^{(2)}_1,\ldots,;\ldots,;
x^{(0)}_{n-1},x^{(1)}_{n-1},x^{(2)}_{n-1},\ldots).
\]
\end{lem}
\proof
For simplicity, we consider the case $n=2$. The isomorphism is such that
the subspace $G_a$ of $G$ corresponds to the space of degree $a$ homogeneous
polynomials in $(x^{(0)}_1,x^{(1)}_1,x^{(2)}_1,\ldots)$.
The isomorphism is given by
\bea
x^{(m_1)}_1\cdots x^{(m_a)}_1\mapsto
S\Bigl((u^{(1)}_1)^{m_1}\cdots(u^{(a)}_1)^{m_a}\Bigr).
\ena
\qed

The basic relations are
\begin{lem}
\bea
&&f_{i\cdots i+l-1}[k_1,\ldots,k_{l-1}]*f_{i+l\cdots i+l+m}
[k_{l+1},\ldots,
k_{l+m}]
\nonumber\\
&&
=f_{i\cdots i+l+m}[k_1,\ldots,
k_{l-1},-1,k_{l+1},\ldots,k_{l+m}],
\nonumber\\\label{QR1}\\
&&f_{i+l\cdots i+l+m}[k_{l+1},\ldots,k_{l+m}]*f_{i\cdots i+l-1}[k_1,\ldots,k_{l-1}]
\nonumber\\
&&
=f_{i\cdots i+l+m}[k_1,\ldots,k_{l-1},
0,k_{l+1},\ldots,k_{l+m}],
\nonumber\\\label{QR2}\\
&&f_i*f_{i\cdots i+m}[k_1,k_2,\ldots,k_m]
=f_{i\cdots i+m}[k_1-1,k_2,\ldots,k_m]*f_i,
\nonumber\\
\label{QR21}\\
&&f_{i+m}*f_{i\cdots i+m}[k_1,\ldots,k_{m-1},k_m-1]=
f_{i\cdots i+m}[k_1,\ldots,k_{m-1},k_m]*f_{i+m},
\nonumber\\
\label{QR6}\\
&&f_{i\cdots i+m}[k_1,\ldots,k_m]*f_{j\cdots j+l}[k'_1,\ldots,k'_l]=
f_{j\cdots j+l}[k'_1,\ldots,k'_l]
\nonumber\\
&&
*f_{i\cdots i+m}[k_1,\ldots,k_m],
\phantom{mmmmmmm}\phantom{mmmmmmm}\hbox{ if $i+m+1<j$.}
\nonumber\\
\label{QR7}
\ena
\end{lem}
The proof is straightforward.

\begin{lem}
\bea
&&f_{i\cdots i+l}[k_1,\ldots,k_l]*
f_{i\cdots i+l+m}[k_1+k_{l+1},\ldots,k_l+k_{l+1},k_{l+1},k_{l+2},\ldots,
\nonumber\\&&
k_{l+m}]
=f_{i\cdots i+l+m}[k_1+k_{l+1},\ldots,k_l+k_{l+1},k_{l+1}-1,k_{l+2},\ldots,k_{l+m}]
\nonumber\\&&
*f_{i\cdots i+l}[k_1,\ldots,k_l],
\nonumber\\
&&\label{REL1}\\
&&f_{i+l\cdots i+l+m}[k_{l+1},\ldots,k_{l+m}]*
f_{i\cdots i+l+m}[k_1,\ldots,k_{l-1},k_l-1,
k_{l+1}+k_l,\ldots,\nonumber\\
&&
\ldots,
k_{l+m}+k_l]
=f_{i\cdots i+l+m}[k_1,\ldots,k_{l-1},k_l,k_{l+1}+k_l,
\ldots,
\nonumber\\
&&\ldots,
k_{l+m}+k_l]*
f_{i+l\cdots i+l+m}[k_{l+1},\ldots,k_{l+m}],
\nonumber\\
&&\label{REL2}\\
&&f_{i\cdots i+l+m}[k_1,\ldots,k_{l-1},k_l,k_{l+1}+k_l,\ldots,k_{l+m}+k_l]
\nonumber\\
&&*f_{i+l\cdots i+l+m+p}
[k_{l+1},\ldots,k_{l+m},-k_l,k_{l+m+2},\ldots,k_{l+m+p}]
\nonumber\\
&&=f_{i+l\cdots i+l+m+p}[k_{l+1},\ldots,k_{l+m},-k_l-1,k_{l+m+2},
\ldots,k_{l+m+p}]
\nonumber\\&&
*f_{i\cdots i+l+m}[k_1,\ldots,k_{l-1},
k_l-1,k_{l+1}+k_l,\ldots,k_{l+m}+k_l].
\nonumber\\
\label{REL3}
\ena
\end{lem}

\noindent
The proof will be given later.

The following are simple consequences of the above.

\begin{prop}
\bea
&&f^{(a)}_{i\cdots i+l}[k_1,\ldots,k_l]*
f^{(b)}_{i\cdots i+l+m}[k_1+k_{l+1},\ldots,k_l+k_{l+1},k_{l+1}+a-1,k_{l+2},\ldots,
\nonumber\\&&\ldots,
k_{l+m}]
=f^{(b)}_{i\cdots i+l+m}[k_1+k_{l+1},\ldots,k_l+k_{l+1},k_{l+1}-1,k_{l+2},\ldots,k_{l+m}]
\nonumber\\&&
*f^{(a)}_{i\cdots i+l}[k_1,\ldots,k_l],
\nonumber\\
\label{PCS1}
\\&&
f^{(a)}_{i+l\cdots i+l+m}[k_{l+1},\ldots,k_{l+m}]*
f^{(b)}_{i\cdots i+l+m}[k_1,\ldots,k_{l-1},k_l-1,k_{l+1}+k_l,\ldots,
\nonumber\\&&
\ldots,
k_{l+m}+k_l]
=f^{(b)}_{i\cdots i+l+m}[k_1,\ldots,k_{l-1},k_l+a-1,k_{l+1}+k_l,\ldots,k_{l+m}+k_l]
\nonumber\\&&
*f^{(a)}_{i+l\cdots i+l+m}[k_{l+1},\ldots,k_{l+m}],\label{PCS2}
\nonumber\\
\\&&f^{(a)}_{i\cdots i+l+m}[k_1,\ldots,k_{l-1},k_l+b-1,k_{l+1}+k_l,\ldots,k_{l+m}+k_l]
\nonumber\\&&
*f^{(b)}_{i+l\cdots i+l+m+p}[k_{l+1},\ldots,k_{l+m},-k_l+a-1,k_{l+m+2},\ldots,k_{l+m+p}]
\nonumber\\&&
=f^{(b)}_{i+l\cdots i+l+m+p}[k_{l+1},\ldots,k_{l+m},-k_l-1,k_{l+m+2},\ldots,k_{l+m+p}]
\nonumber\\&&
*f^{(a)}_{i\cdots i+l+m}[k_1,\ldots,k_{l-1},k_l-1,k_{l+1}+k_l,\ldots,k_{l+m}+k_l].
\nonumber\\
\label{FR3}
\ena
\end{prop}

The following is also valid. However, the general case for $a,b$ does not
follow from the special case $a=b=1$.
\begin{prop}
\bea
&&f^{(a)}_{i+l\cdots i+l+m}[k_{l+1},\ldots,k_{l+m}]
*f^{(b)}_{i\cdots i+l+m+p}[k_1,\ldots,\nonumber\\&&
\ldots,
k_{l-1},k_l-1,
k_{l+1}+k_l,\ldots,k_{l+m}+k_l,k_l+a-1,k_{l+m+2},\ldots,k_{l+m+p}]
\nonumber\\&&
=f^{(b)}_{i\cdots i+l+m+p}[k_1,\ldots,k_{l-1},k_l+a-1,
k_{l+1}+k_l,\ldots,k_{l+m}+k_l,
\nonumber\\&&
k_l-1,k_{l+m+2},\ldots,k_{l+m+p}]
*f^{(a)}_{i+l\cdots i+l+m}[k_{l+1},\ldots,k_{l+m}]
\nonumber\\\label{FR4}
\ena
\end{prop}

We use these relations for integer parameters $k_i$.
However, they are valid without this restriction because
the general case follows from the integer case.

Let us derive (\ref{REL1}). The other cases follow similarly.
Without loss of generality one can assume that $m=1$.
Using (\ref{COM}) and (\ref{QR2}) we have
\bea
&&
f_{i+l}*f_{i\cdots i+l-1}[k_1,\ldots,k_{l-1}]
*f_{i\cdots i+l}[k_1+k_l,\ldots,k_{l-1}+k_l,k_l]\nonumber\\&&
=f_{i\cdots i+l}[k_1+k_l,\ldots,k_{l-1}+k_l,k_l]
*f_{i+l}*f_{i\cdots i+l-1}[k_1,\ldots,k_{l-1}].
\nonumber\\
\ena
Using (\ref{QR6}) we have
\bea
&&
f_{i+l}*\Bigl(f_{i\cdots i+l-1}[k_1,\ldots,k_{l-1}]*
f_{i\cdots i+l}[k_1+k_l,\ldots,k_{l-1}+k_l,k_l]\nonumber\\&&
-f_{i\cdots i+l}[k_1+k_l,\ldots,k_{l-1}+k_l,k_l-1]
*f_{i\cdots i+l-1}[k_1,\ldots,k_{l-1}]\Bigr)=0.
\nonumber\\
\ena
Since $f_{i+l}$ is not a zero divisor, we have (\ref{REL1}).

Combining all these, in particular
(\ref{PCS1}) and (\ref{PCS2}), 
we arrive at the formulas which we need for the quadratic relations 
of the screening operators. 
\begin{prop}
\bea
&&f^{(a+b)}_{i\cdots i+l}[k_1,\ldots,k_l]*
f^{(b)}_{i+l+1\cdots i+l+m}[k_{l+2},\ldots,k_{l+m}]
\nonumber\\&&
=f^{(b)}_{i\cdots i+l+m}[k_1-a,\ldots,k_l-a,-a-1,k_{l+2},\ldots,k_{l+m}]*
f^{(a)}_{i\cdots i+l}[k_1,\ldots,k_l],
\nonumber\\
\label{CS1}
\\&&
f^{(a)}_{i\cdots i+l}[k_1,\ldots,k_l]*
f^{(a+b)}_{i+l+1\cdots i+l+m}[k_{l+2},\ldots,k_{l+m}]
\nonumber\\&&
=f^{(b)}_{i+l+1\cdots i+l+m}[k_{l+2},\ldots,k_{l+m}]*
f^{(a)}_{i\cdots i+l+m}[k_1,\ldots,k_l,-b-1,k_{l+2}-b,\ldots,
\nonumber\\&&
\ldots,k_{l+m}-b],
\nonumber\\
\label{CS2}
\\&&
f^{(a)}_{i+l+1\cdots i+l+m}[k_{l+2},\ldots,k_{l+m}]*
f^{(a+b)}_{i\cdots i+l}[k_1,\ldots,k_l]
\nonumber\\&&
=f^{(b)}_{i\cdots i+l}[k_1-a,\ldots,k_l-a]*
f^{(a)}_{i\cdots i+l+m}[k_1,\ldots,k_l,a+b-1,k_{l+2},\ldots,k_{l+m}],
\nonumber\\
\label{CS3}
\\&&
f^{(a+b)}_{i+l+1\cdots i+l+m}[k_{l+2},\ldots,k_{l+m}]*
f^{(b)}_{i\cdots i+l}[k_1,\ldots,k_l]
\nonumber\\&&
=f^{(b)}_{i\cdots i+l+m}[k_1,\ldots,k_l,a+b-1,k_{l+2},\ldots,k_{l+m}]
*f^{(a)}_{i+l+1\cdots i+l+m}[k_{l+2}-b,\ldots,
\nonumber\\&&
\ldots,k_{l+m}-b].
\nonumber\\\label{CS4}
\ena
\end{prop}
\proof
Let us derive (\ref{CS1}). The other cases can be proven similarly.
Using (\ref{MSC}), (\ref{QR1}), (\ref{PCS1}), we have
\bea
&&f^{(a+b)}_{i\cdots i+l}[k_1,\ldots,k_l]*
f^{(b)}_{i+l+1\cdots i+l+m}[k_{l+2},\ldots,k_{l+m}]
\nonumber\\&&
=f^{(a+b-1)}_{i\cdots i+l}[k_1,\ldots,k_l]*
f_{i\cdots i+l}[k_1-a-b+1,\ldots,k_l-a-b+1]
\nonumber\\&&
*f_{i+l+1\cdots i+l+m}[k_{l+2}-b+1,\ldots,k_{l+m}-b+1]*
f^{(b-1)}_{i+l+1\cdots i+l+m}[k_{l+2},\ldots,k_{l+m}]
\nonumber\\&&
=f^{(a+b-1)}_{i\cdots i+l}[k_1,\ldots,k_l]
*f_{i\cdots i+l+m}[k_1-a-b+1,\ldots,
k_l-a-b+1,-1,
\nonumber\\&&
k_{l+2}-b+1,\ldots,k_{l+m}-b+1]
*f^{(b-1)}_{i+l+1\cdots i+l+m}[k_{l+2},\ldots,k_{l+m}]
\nonumber\\&&
=f_{i\cdots i+l+m}[k_1-a-b+1,\ldots,k_l-a-b+1,-a-b,k_{l+2}-b+1,\ldots,
\nonumber\\&&
\ldots k_{l+m}-b+1]
*f^{(a+b-1)}_{i\cdots i+l}[k_1,\ldots,k_l]*
f^{(b-1)}_{i+l+1\cdots i+l+m}[k_{l+2},\ldots,k_{l+m}]
\nonumber\\&&
=f^{(b)}_{i\cdots i+l+m}[k_1-a,\ldots,k_l-a,-a-1,k_{l+2},\ldots,k_{l+m}]
*f^{(a)}_{i\cdots i+l}[k_1,\ldots,k_l].
\nonumber\\
\ena
\qed

\setcounter{equation}{0}
\section{Commutativity with DWA.}\label{app:3}
   
In this section we prove the Lemma \ref{LEMDWA2}. 
It would be more convenient for us to use here the 
"multiplicative" variable $z$ instead of $u$. 
For this reason, let us define the theta function 
\bea
[[z]]\equiv[u]\quad ,\quad z=x^{2u} \
\ena
having the periodicity property $[[zx^{2r}]]=-[[z]]$.
Abusing the notations, let us use the same symbol 
$f_\alpha^{(a)}(z_i^{(1)},\cdots,z_{i+m}^{(a)})$  
for the function 
$f_\alpha^{(a)}(u_i^{(1)},\cdots,u_{i+m}^{(a)})$ where
$z_{j}^{(b)}=x^{2u_{j}^{(b)}}$.
The screening operator 
$\bX_\alpha(\lambda)$ in the notations 
\eqref{SET1}-\eqref{SET4} is given by \eqref{eqn:scr}, i.e., in the
multiplicative variables,
\bea
&&\overline X_\alpha(\lambda)=
\oint\cdots\oint\prod_{1\leq b\leq a\atop i\leq j\leq i+m}
{dz^{(b)}_{j}\over2\pi iz^{(b)}_{j}}\nonumber\\
&&\times 
\Biggl({U_i(z^{(1)}_i)\over[[z^{(1)}_i/x]]}
\cdots{U_i(z^{(a)}_i)\over[[z^{(a)}_i/x]]}\Biggr)
\cdots
\Biggl({U_{i+m}(z^{(1)}_{i+m})\over[[z^{(1)}_{i+m}/x]]}
\cdots {U_{i+m}(z^{(a)}_{i+m})\over[[z^{(a)}_{i+m}/x]]}\Biggr)
\nonumber\\
&&
\times\prod_{1\leq b<c\leq a\atop i\leq j\leq i+m}
{[[z^{(b)}_{j}/z^{(c)}_{j}]]
\over[[z^{(b)}_{j}/x^2z^{(c)}_{j}]]}
\prod_{1\leq b,c\leq a\atop i\leq j\leq i+m-1}
{[[z^{(b)}_{j}/z^{(c)}_{j+1}]]
\over[[z^{(b)}_{j}x/z^{(c)}_{j+1}]]}
\nonumber\\
&&\times f^{(a)}_\alpha(z^{(1)}_i,\ldots,z^{(a)}_i,\ldots,
z^{(1)}_{i+m},\ldots,z^{(a)}_{i+m}).
\label{dwaSCR}\ena

The non-trivial couplings between 
the screening currents with each 
others and with the field $\Lambda_j(z)$ are:
\bea
&&U_j(z)U_{j+1}(w)=
z^{-\frac{r-1}{r}}s(w/z):U_{j}(z)U_{j+1}(w):\ ,\cr
&&U_{j+1}(z)U_j(w)=
z^{-\frac{r-1}{r}}s(w/z):U_{j+1}(z)U_j(w):\ ,\cr
&&
U_j(z)U_{j}(w)=z^{2\frac{r-1}{r}}t(w/z):U_{j}(z)U_{j}(w):\ ,\cr
&&
\Lambda_j(z)U_j(w)=x^{-2(r-1)}
\tau_j(w/z):\Lambda_j(z)U_j(w):\ , \cr
&&
U_j(w)\Lambda_j(z)=x^{-2(r-1)}
\tau_j(w/z):\Lambda_j(z)U_j(w):\ , \cr
&&
U_j(z)\Lambda_{j+1}(w)=
\hat{\tau}_{j}(w/z):\Lambda_{j+1}(w)U_j(z):\ , \cr
&&
\Lambda_{j+1}(w)U_j(z)=
\hat{\tau}_{j}(w/z):\Lambda_{j+1}(w)U_j(z):\ , \label{C.4}
\ena
where 
\bea
&&s(z)=\frac{(x^{2r-1}z;x^{2r})_\infty}{(xz;x^{2r})_\infty}\ ,\cr
&&t(z)=(1-z)\frac{(x^{2}z;x^{2r})_\infty}{(x^{2r-2}z;x^{2r})_\infty}
\ , \cr
&&\tau_j(z)=\frac{1-zx^{r+j-2}}{1-zx^{j-r}}\ , \cr
&&\hat{\tau}_{j}(z)=\frac{1-zx^{r-j-2}}{1-zx^{-r-j}}\ . 
\label{cop}\ena

To be precise, the $\tau_j(w/z)$ (resp.
$\hat{\tau}_{j}(w/z)$)
in the formula for
$U_j(w)\Lambda_j(z)$ (resp.
$\Lambda_{j+1}(w)U_j(z)$)
must be understood as the expansion
in $z/w$.

The screening currents commute with the
generator ${\cal W}^{(1)}_t$ up to a total difference \cite{qWN,FeFr95}.
\begin{lem}\label{LEMdwa1}
\begin{equation}
[{\cal{W}}^{(1)}_t, U_j(z)]
=(x^{-2r+2}-1)(zx^{j-r})^{t}:\Lambda_j(zx^{j-r})U_j(z): 
-\{z\rightarrow zx^{2r} \}\ . 
\end{equation}
\end{lem}
According to this Lemma we have 
\bea
&&(x^{-2r+2}-1)^{-1}[{\cal W}^{(1)}_t,U_i(z_{i}^{(1)})
\cdots U_{i+m}(z_{i+m}^{(a)})]
\nonumber\\
&&
=\sum_{i\leq j\leq i+m\atop 1\leq f\leq a} 
\Bigl((z_{j}^{(f)}x^{j-r})^{t}
U_{i}(z_{i}^{(1)})\cdots 
:\Lambda_{j}(z_{j}^{(f)}x^{j-r})U_{j}(z_{j}^{(f)}):
\cdots U_{i+m}(z_{i+m}^{(a)}) 
\nonumber\\
&&-\{z_{j}^{(f)}\rightarrow x^{2r}z_{j}^{(f)} \} \Bigr).\label{STP1}
\ena
By normal-ordering the operator part we obtain an expression 
involving various functions and operators. 
For our purpose it is useful to introduce the following objects. 
Let $0\le s\le m$, $1 \le f_0,\cdots,f_s \le a$ and $i\le k\le i+m-s$. 
(The case of our immediate interest above is $s=0$. However, the general
case will be necessary as we proceed.)
We define
\bea
&&I^{(f_0,\ldots,f_s)}_{k,\ldots,k+s}\nonumber\\
&=&
\left(\oint_{|w|=1}-\oint_{|w|=x^{2r-s}}\right)
\iv{w}
J^{(f_0)}_k(z^{(1)}_i,\ldots,z^{(a)}_{i+m})
\Bigg|_\rst
\ena
where
\bea
&&J^{(f_0)}_k(\zvar)=U^{(f_0)}_k(\zvar)\nonumber\\
&&\phantom{mmm}\times F(\zvar)T^{(f_0)}_k(\zvar)S^{(f_0)}_k(\zvar)\nonumber
\ena
and $U^{(f_0)}_k$, $F$, $T^{(f_0)}_k$ and $S^{(f_0)}_k$
are given below.
\bea
U^{(f_0)}_k(\zvar)
&=&(z^{(f_0)}_kx^{k-r})^t
:\Lambda_k(z^{(f_0)}_kx^{k-r})
\prod_{1\le b\le a\atop i\le j\le i+m}
U_j(z^{(b)}_j):,
\nonumber\\
&&\\
F(\zvar)&=&
\prd<[[\z bj/\z cj]]\prd,[[\z bj/\z c{j+1}]]\nonumber\\
&&\times
f^{(a)}_\alpha(\zvar)
\prod_{1\le b\le a\atop i\le j\le i+m}{1\over[[\z bj/x]]},\nonumber\\
T^{(f_0)}_k(\zvar)&=&
\prd<{(\z bj)^{{2(r-1)\over r}}t(\z cj/\z bj)\over[[\z bj/x^2\z cj]]}
\nonumber\\
&&\times
\prod_{1\le b\le a\atop b\not=f_0}
x^{-2(r-1)}\tau_k(\z bkx^{r-k}/\z {f_0}k),\\
S^{(f_0)}_k(\zvar)&=&
\prod_{1\le b,c\le a\atop i\le j\le i+m-1}}
{(\z bj)^{-{r-1\over r}}s(\z c{j+1}/\z bj)
\over[[x\z bj/\z c{j+1}]]\nonumber\\
&&\times\prod_{1\le b\le a}
\hat{\tau}_{k-1}(\z {f_0}k/x^{r-k}\z b{k-1}).
\ena

The restriction $\rst$ is regular except for
the functions
$s(\z {f_{p+1}}{k+p+1}/\z {f_p}{k+p})$.
For these singular terms we use the convention
\bea\label{THREESTARS}
s(z^{-1})\Bigg|_{z=x}&{\buildrel{\rm def}\over=}&
{\rm Res\,}_{z=x}s(z^{-1})\iv z=
{(x^{2(r-1)};x^{2r})_\infty\over(x^{2r};x^{2r})_\infty}.
\ena

Set
\bea
I(s)&=&
\sum_{k=i}^{i+m}\sum_{f_0,\ldots,f_s=1}^a\oint\cdots\oint
\prod_{{\scriptstyle1\le b\le a\atop i\le j \le i+m}
\atop(b,j)\not=(f_0,k),\ldots,(f_s,k+s)}
\iv{\z bj}
I^{(f_0,\ldots,f_s)}_{k,\ldots,k+s}.
\ena
The contour for $\z bj$ is $|\z bj|=1$.

The induction goes as follows.

\noindent
Step 1\quad $[{\cal W}^{(1)}_t,{\overline X}_\alpha(\lambda)]=I(0)$

\noindent
Step 2\quad $I(s)=I(s+1)$ $(0\le s\le m-1)$

\noindent
Step 3\quad $I(m)=0$

Step 1 follows from (\ref{dwaSCR}), (\ref{C.4}) and (\ref{STP1})

We will show that the only poles of $J^{(f_0)}_k(\zvar)\Bigg |_\rst$
between the contours $|w|=1$ and $|w|=x^{2r-s}$
are
\bea
&(A)&\quad w=x\z b{k+s+1}\quad(i\le k\le i+m-s-1,1\le b\le a),\nonumber\\
&(B)&\quad w=x^{2r-s-1}\z b{k-1}\quad(i+1\le k\le i+m-s,1\le b\le a).
\label{TWOSTARS}
\ena
In particular, Step 3 follows.

Let us abbreviate $U^{(f_0)}_k(\zvar)\Bigg|_\rst$, etc.,
by $U^{(f_0)}_k\Bigg|_s$.
\begin{itemize}
\item $\ab U$ has no poles. This is because the term is normal-ordered.
\item $F\Bigg|_s$ has no poles. This is proved in Proposition \ref{prop:scr}.
\item
$T^{(f_0)}_k\Bigg |_s$ has no poles in the region
\bea\label{ONESTAR}
x^{2r-s}<|w|<1.
\ena

\noindent(a) Consider the poles of 
$[[\z bj/x^2\z cj]]^{-1}\Bigg |_{\z bj=x^pw}$.
They are situated at $w=$
\newline
$x^{2r\nu+2-p}\z cj$ $(\nu\in\Z)$. However
the factor $t(\z cj/\z bj)\Bigg|_{\z bj=x^pw}$
cancels the poles for $\nu\ge0$. Therefore, there are no such poles in the
region (\ref{ONESTAR}).

\noindent(b) Consider the poles of
$[[\z bj/x^2\z cj]]^{-1}\Bigg|_{\z cj=x^pw}$.
The possible poles in (\ref{ONESTAR}) are
$w=x^{2r-s-2}\z bk$ $(b<f_0)$ or $w=x^{2r-s-1}\z b{k+1}$.
The former is canceled by $\tau_k(\z bkx^{r-k}/x^sw)$,
and the latter is canceled by the property (P5) of $f^{(a)}_\alpha$.

\noindent(c) Consider the poles of $t(\z cj/\z bj)\Bigg|_{\z bj=x^pw}$.
The possible poles in (\ref{ONESTAR}) are
$w=x^{2r-s-1}\z c{k+1}$ or $w=x^{2r-s-2}\z c{k+1}]$ $(f_0<c)$.
The former is canceled by the property (P5) of $f^{(a)}_\alpha$
and the latter is canceled by $\tau_k(\z ckx^{r-k}/x^sw)$.

\noindent(d) There is no poles of $t(\z cj/\z bj)\Bigg|_{\z cj=x^pw}$
in (\ref{ONESTAR}).

\noindent(e) There is no pole of $\tau_k(\z bkx^{r-k}/x^sw)$ in (\ref{ONESTAR}).
\item $S^{(f_0)}_k\Bigg|_s$ has two types of poles (A)
and (B) of (\ref{TWOSTARS}) arising from the factors
$s(\z c{k+s+1}/\z {f_s}{k+s})\Bigg|_{\z{f_s}{k+s}=w}$
and $\hat{\tau}_{k-1}(\z{f_0}k/x^{r-k}\z b{k-1})\Bigg|_{\z {f_0}k=x^sw}$,
respectively.

\noindent(a) Consider the poles of $s(\z c{j+1}/\z bj)\Bigg|_{\z bj=x^pw}$.
They are situated at $w=x^{p-s+1+2r\nu}\z c{k+p+1}$ $(\nu\in\Z_{\ge0})$.
The only  poles in (\ref{ONESTAR}) are $w=x\z c{k+s+1}$ $(1\le c\le a)$.

\noindent(b) The pole of
$\hat{\tau}_{k-1}(\z {f_0}k/x^{r-k}\z b{k-1})\Bigg|_{\z{f_0}k=x^sw}$
is at $w=x^{2r-s-1}\z b{k-1}$.

\noindent(c) Consider the poles of
$[[x\z bj/\z c{j+1}]]^{-1}\Bigg|_{\z bj=x^{s-p}w}$. They are 
situated at
$w=x^{2r+p-s-1}\z c{k+p}$ $(\nu\in\Z)$.
The only pole in (\ref{ONESTAR}) is $w=x^{2r-s-1}\z ck$.
This is canceled by the property (P5) of $f^{(a)}_\alpha$.

\noindent(d) Consider the poles of
$[[x\z bj/\z c{j+1}]]^{-1}\Bigg|_{\z c{j+1}=z^{s-p}w}$. They are 
situated at
$w=x^{2r\nu+p-s+1}\z b{k+p-1}$ $(\nu\in\Z)$. The only pole in
(\ref{ONESTAR}) is $w=x\z b{k+s-1}$. This is canceled by the property (P5)
of $f^{(a)}_\alpha$.
\end{itemize}

\noindent
We finished checking (\ref{TWOSTARS}).

We will conclude the induction by showing Step 2.
We take the residues of $J^{(f_0)}_k(\zvar)\Bigg|_{\rst}$
at (A) and (B) of (\ref{TWOSTARS}).

\noindent(A):
We take the residue at $w=x\z {f_{s+1}}{k+s+1}$,
and then rename $\z{f_{s+1}}{k+s+1}$ to $w$ in order to compare it with
$J^{(f_0)}_k(\zvar)\Bigg|_{z^{(f_p)}_{k+p}=x^{s+1-p}w\phantom{m}(0\le p\le s+1)}$.
Except for the factor $s(\z{f_{s+1}}{k+s+1}/\z{f_s}{k+s})$,
this procedure is equivalent to changing the substitution rule
$\z{f_p}{k+p}=x^{s-p}w$ $(0\le p\le s)$
to
$\z{f_p}{k+p}=x^{s+1-p}w$ $(0\le p\le s+1)$.
The residue of $s(\z{f_{s+1}}{k+s+1}/w)\iv w$ at $w=x\z{f_{s+1}}{k+s+1}$
gives $s(x^{-1})$ in the convention of (\ref{THREESTARS}).
Thus we get the term in $I^{f_0,\ldots,f_{s+1}}_{k,\ldots,k+s+1}$
that is corresponding to the cycle $|w|=1$.

\noindent(B): It is convenient to consider the residue of
$J^{(f_1)}_{k+1}(\zvar)_{z^{(f_{p+1})}_{k+1+p}
=x^{s-p}w\phantom{m}}$, $(0\le p\le s)$
at $w=x^{2r-s-1}\z{f_0}k$.
We take the residue and then rename
$\z{f_0}k$ to $x^{-2r+s+1}w$. We will compare

\noindent
(I)\quad
$\displaystyle J^{(f_1)}_{k+1}(\zvar)
\Bigg|_{z^{(f_{p+1})}_{k+1+p}=x^{s-p}w\phantom{m}
(0\le p\le s),\phantom{m}\z{f_0}k=x^{-2r+s+1}w}$

\noindent
and

\noindent
(II)\quad
$\displaystyle J^{(f_0)}_k(\zvar)
\Bigg|_{z^{(f_p)}_{k+p}=x^{s+1-p}w\phantom{m}(0\le p\le s+1)}$

\noindent
taking care of the pole of
$\hat\tau_k(\z{f_1}{k+1}/x^{r-k-1}\z{f_0}k)\Bigg|_{\z{f_1}{k+1}=x^sw}$
at $w=x^{2r-s-1\z{f_0}k}$.

Let us abbreviate the restrictions (I) and (II) by $\Bigg|_I$
and $\Bigg|_{II}$, respectively.
We will compute the ratios of the corresponding terms in (I) and (II).
\begin{itemize}
\item (U): $U^{(f_1)}_{k+1}\Bigg|_I/U^{(f_0)}_k\Bigg|_{II}=x^{-2(r-1)}$

This follows from the identity
\[
x^{2(r-1)}:\Lambda_{k+1}(x^{r+k}z)U_k(z):=
:\Lambda_k(x^{r+k}z)U_k(x^{2r}z):.
\]
\end{itemize}

Set
\[
T_0=\prd<[[\z bj/x^2\z cj]]
\]
and
\[
S_0=\prd,[[x\z bj/\z c{j+1}]].
\]
\begin{itemize}
\item (F): $F(T_0S_0)^{-1}\Bigg|_I/F(T_0S_0)^{-1}\Bigg|_{II}=1$

The signs arising from the periodicity
$[[x^{2r}z]]=-[[z]]$ cancel out.

\item (T): $T^{(f_1)}_{k+1}T_0\Bigg|_I/T^{(f_0)}_kT_0\Bigg|_{II}=$
\[
\prod_{b\not=f_1}
{1-x^{2r-s-2}\z b{k+1}/w\over1-x^{-s}\z b{k+1}/w}
\prod_{b\not=f_0}x^{-4(r-1)}
{1-x^{2r-s-1}\z bk/w\over1-x^{-s+1}\z bk/w}
\]

We used 
\[
{t(x^{2r}z)\over t(z)}={(1-x^{2(r-1)}z)(1-x^{2r}z)\over{(1-z)(1-x^2z)}}.
\]

\item (S1):
\[
{\rm Res}_{w=x^{2r-s-1}\z{f_0}k}
\hat\tau_k(x^{s+k+1-r}w/\z{f_0}k)\iv w
=x^{2(r-1)}-1.
\]

\item (S2):
\bea
&&{\prod_{b\not=f_0}\hat\tau_k(\z{f_1}{k+1}/x^{r-k-1}\z bk)\Bigg|_{\z{f_1}{k+1}=x^sw}
\over
\prod_b\hat\tau_{k-1}(\z{f_0}k/x^{r-k}\z b{k-1})\Bigg|_{\z{f_0}k=x^{s+1}w}}
\nonumber\\
&&=
\prod_{b\not=f_0}x^{2(r-1)}
{1-x^{-s+1}\z bk/w\over1-x^{2r-s-1}\z bk/w}
\prod_b{1-x^{-2r+s+2}w/\z b{k-1}\over1-x^sw/\z b{k-1}}
\ena
\end{itemize}

In the restriction (I) we must pay a special attention to the factor
$(\z{f_0}k)^{-{r-1\over r}}$
\newline
$\times s(\z{f_1}{k+1}/\z{f_0}k)$. We have
\bea
(\z{f_0}k)^{-{r-1\over r}}s(\z{f_1}{k+1}/\z{f_0}k)
\Bigg|_{\z{f_1}{k+1}=x^sw\atop\z{f_0}k=x^{-2r+s+1}w}
&=&
x^{2(r-1)}
(x^{s+1}w)^{-{r-1\over r}}
s(x^{2r-1})\nonumber\\
&=&
{x^{2(r-1)}(x^{s+1}w)^{-{r-1\over r}}\over1-x^{2(r-1)}}s(x^{-1})\nonumber
\ena
where
 $s(x^{-1})$ is  in the sense of (\ref{THREESTARS}).
Taking this into account and collecting (S1) and (S2)
we have
\begin{itemize}
\item (S): $S^{(f_1)}_{k+1}S_0\Bigg|_I/S^{(f_0)}_kS_0\Bigg|_{II}=$
\bea
&&-x^{2(r-1)}
{1-x^{-s+1}\z bk/w\over1-x^{2r-s-1}\z bk/w}
\prod_{b\not=f_1}x^{2(r-1)}{1-x^{-s}\z b{k+1}/w\over1-x^{2r-s-2}\z b{k+1}/w}.
\nonumber
\ena
\end{itemize}

{}From (U), (F), (T) and (S) we can conclude that the residues of
\[
J^{(f_1)}_{k+1}(\zvar)\Bigg|_{z^{(f_{p+1})}_{k+p+1}=x^{s-p}w}
\phantom{m}(0\le p\le s)
\]
at $w=x^{2r-s-1}\z{f_0}k$ gives the term in
$I^{(f_0,\ldots,f_{s+1})}_{k,\ldots,k+s+1}$
that is corresponding to the cycle $|w|=x^{2r-s-1}$.
This completes the proof of the commutativity
$[{\cal W}^{(1)}_t,{\overline X}_\alpha(\lambda)]=0$.

\end{document}